\documentclass[10pt]{article}

\usepackage{graphicx,hyperref}
\usepackage{amsmath}
\usepackage[version=4]{mhchem}
\usepackage{siunitx}
\usepackage{longtable,tabularx}
\usepackage{pstricks}
\usepackage{amsfonts}
\usepackage{latexsym}
\usepackage{enumerate}
\usepackage{epsfig}
\usepackage{subfig}
\usepackage{framed,fancybox}
\usepackage{bbm}
\usepackage{multirow}
\usepackage[T1]{fontenc}
\usepackage{threeparttable}
\usepackage{rotating}
\usepackage{amssymb,amscd}
\makeatletter
\def\myVCENTER#1{\vcenter{\hbox{$\m@th#1$}}}
\makeatother

\long\def\symbolfootnote[#1]#2{\begingroup\def\thefootnote{\fnsymbol{footnote}}\footnote[#1]{#2}\endgroup}

\definecolor{shadecolor}{gray}{0.99}
{\endMakeFramed}

\newenvironment{shadedframe}{%
 \MakeFramed {\FrameRestore}}
{\endMakeFramed}

\definecolor{shadecolor}{gray}{0.99}
\long\def\symbolfootnote[#1]#2{\begingroup\def\thefootnote{\fnsymbol{footnote}}\footnote[#1]{#2}\endgroup}
\def\qed{\hfill{$\vcenter{\hrule height1pt \hbox{\vrule width1pt height5pt
    \kern5pt \vrule width1pt} \hrule height1pt}$} \medskip}
\newcommand{\m}[1]{{\bf{#1}}}
\newcommand{\g}[1]{\boldsymbol #1}
\newcommand{\bb}[1]{\mathbb #1}
\newcommand{\tr}{^{\sf T}}
\newcommand{\C}[1]{{\cal {#1}}}

\title{\bf Mesh Refinement Method for Solving Bang-Bang Optimal Control Problems Using Direct Collocation}

\author{Yunus M. Agamawi\thanks{Ph.D.~Student, Department of Mechanical and Aerospace Engineering.  E-mail:  yagamawi@ufl.edu.}\\William W.~Hager\thanks{Distinguished Professor, Department of Mathematics.  E-mail:  hager@ufl.edu.}\\ Anil V. Rao\thanks{Associate Professor, Department of Mechanical and Aerospace Engineering, Erich Farber Faculty Fellow and University Term Professor.  E-mail: anilvrao@ufl.edu.  Associate Fellow, AIAA.  Corresponding Author.} \vspace{12pt} \\ {\em University of  Florida, Gainesville, FL, 32611}}

\date{}
\begin{document}

\maketitle
\renewcommand{\baselinestretch}{1}
\normalsize\normalfont 
\begin{abstract}
  \noindent  A mesh refinement method is developed for solving bang-bang optimal control problems using direct collocation.  The method starts by finding a solution on a coarse mesh.  Using this initial solution, the method then determines automatically if the Hamiltonian is linear with respect to the control, and, if so, estimates the locations of the discontinuities in the control.  The switch times are estimated by determining the roots of the switching functions, where the switching functions are determined using estimates of the state and costate obtained from the collocation method.   The accuracy of the switch times is then improved on subsequent meshes by dividing the original optimal control problem into multiple domains and including variables that define the locations of the switch times.  While in principle any collocation method can be used, in this research the previously developed Legendre-Gauss-Radau collocation method is employed because it provides an accurate approximation of the costate which in turn improves the approximation of the switching functions.  The method of this paper is designed to be used with a previously developed mesh refinement method in order to accurately approximate the solution in segments where the solution is smooth.  The method is demonstrated on three examples where it is shown to accurately determine the switching structure of a bang-bang optimal control problem.  When compared with previously developed mesh refinement methods, the results demonstrate that the method developed in this paper improves computational efficiency when solving bang-bang optimal control problems.  
\end{abstract}

\section*{Nomenclature}


{\renewcommand\arraystretch{1.0}
\noindent\begin{longtable*}{lcl}
$\m{a}$	& = & vector field for right-hand side of dynamics\\
$\m{A}$	& = & matrix defining vector field for right-hand side of dynamics at collocation points \\
$\m{b}$	& = & vector field for boundary conditions \\
$\m{c}$	& = & vector field for path constraints	\\
$\m{C}$	& = & matrix defining vector field for path constraints at collocation points \\
$\m{D}$	& = & Legendre-Gauss-Radau differentiation matrix \\
$\m{D}_{N+1}$	& = & last column of Legendre-Gauss-Radau differentiation matrix \\
$f$	& = & $\m{v}$ independent contribution of $\C{H}$	\\
$\C{H}$	& = & Hamiltonian of optimal control problem	\\
$\C{J}$	& = & objective functional \\
$K$		& = & number of mesh intervals used	\\
$\ell_j^{(k)}$ & = & Lagrange polynomial $j$ of mesh interval $k$ \\
$\m{\C{L}}$	& = & vector field for integrand appearing in Lagrange cost \\
$\m{L}$	& = & matrix defining vector field for integrand at collocation points \\
$\C{M}$	& = & Mayer cost \\
$M$			& = & number of mesh iterations	\\
$n_y$	& = & number of state components	\\
$n_u$	& = & number of control components 	\\
$n_c$	& = & number of path constraints	\\
$n_s$	& = & number of switching time parameters	\\
$N_k$	& = & number of collocation points used in mesh interval $k$	\\
$N$		& = & total number of collocation points used	\\
$N_f$	& = & total number of collocation points used for problem on final mesh	\\
$\C{P}_d$	& = & $d^{th}$ time domain	\\
$Q$	& = & number of constant control domains used to approximate bang-bang control profile	\\
$\C{S}_k$ & = & mesh interval $k$ \\
$t_0$	& = & initial time \\
$t_f$ 	& = & final time \\
$t_s^{[S]}$ & = & $S^{th}$ switch time parameter	\\
$t_d^{i}$	& = & estimated discontinuity time for $i^{th}$ control component	\\
$t_{\sigma}^i$ & = & midpoint time of two adjacent points containing a change in sign of $\sigma_i$ \\
$t_{u}^i$ & = & midpoint time of two adjacent points containing largest absolute difference \\
 & & in value of $i^{th}$ control component \\
$\C{T}$		& = & computation time	\\
$\m{u}(t)$	& = & control in time horizon \\
$\m{U}$	& = & matrix of control parameterization at collocation points \\
$\m{v}$	& = & linearly dependent components of control in Hamiltonian	\\
$\m{w}$	& = & vector of corresponding LGR weights at collocation points \\
$\m{W}$	& = & diagonal matrix of corresponding LGR weights at collocation points \\
$\m{y}(t)$	& = & state in time horizon \\
$\m{Y}(\tau)$ & = & state approximation in $\tau$ domain 	\\
$\m{Y}$	& = & matrix of state approximation at discretized points \\
$\m{z}$	& = & nonlinearly dependent components of control in Hamiltonian	\\
$\g{\Delta}$	& = & defect constraints matrix \\
$\g{\lambda}(t)$	& = & costate in time horizon \\
$\g{\lambda}$	& = & matrix of costate estimates at collocation points \\
$\g{\lambda}_{N+1}$	& = & vector of costate estimates at non-collocated end point \\
$\g{\Lambda}$	& = & matrix of NLP multipliers corresponding to defect constraints at collocation points \\
$\g{\sigma}$	& = & vector field for switching functions of $\C{H}$	\\
$\tau$ & = & domain used for Legendre-Gauss-Radau collocation \\
$\tau_j^{(k)}$ & = & support point $j$ of mesh interval $k$ \\
\end{longtable*}}

\renewcommand{\baselinestretch}{1}
\normalsize\normalfont 

\section{Introduction \label{sect:intro}}

Optimal control problems arise frequently in many engineering applications due to the need to optimize performance of a controlled dynamical system.  In general, optimal control problems do not have analytic solutions and, thus, must be solved numerically.  Numerical methods for optimal control fall into two broad categories: indirect methods and direct methods.  In an indirect method, the first-order variational optimality conditions are derived, and the optimal control problem is converted to a Hamiltonian boundary-value problem (HBVP).  The HBVP is then solved numerically using a differential-algebraic equation solver.  In a direct method, the state and control are approximated, and the optimal control problem is transcribed into a finite-dimensional nonlinear programming problem (NLP) \cite{Betts3}.  The NLP is then solved numerically using well-developed software \cite{Gill1,Gill2,Biegler2}.  

Over the past two decades, a particular class of direct methods, called direct collocation methods, has been used extensively for solving continuous optimal control problems.  A direct collocation method is an implicit simulation method where the state and control are parameterized, and the constraints in the continuous optimal control problem are enforced at a specially chosen set of collocation points.  Traditional direct collocation methods take the form of an $h$ method (for example, Euler or Runge-Kutta methods) where the domain of interest is divided into a mesh, the state is approximated using the same fixed-degree polynomial in each mesh interval, and convergence is achieved by increasing the number and placement of the mesh points \cite{Betts3}.   In contrast to an $h$ method, in recent years so-called $p$ methods have been developed.  In a $p$ method, the number of intervals is fixed, and convergence is achieved by increasing the degree of the approximation in each interval.  To achieve maximum effectiveness, $p$ methods have been developed using {\em orthogonal collocation at Gaussian quadrature points} \cite{Reddien1,Cuthrell2,Elnagar1,Benson2,Rao8,Gong3,Kameswaran1,Garg1,Garg2,Garg3,Darby2,Darby3,Francolin2014a,Patterson2014,Patterson2015,HagerHouRao15a,HagerHouRao16a,HagerLiuMohapatraWangRao18a,DuChenHager2019,HagerHouMohapatraRaoWang2019}.   For problems whose solutions are smooth and well-behaved, a Gaussian quadrature orthogonal collocation method converges at an exponential rate \cite{HagerHouRao15a,HagerHouRao16a,HagerLiuMohapatraWangRao18a,DuChenHager2019,HagerHouMohapatraRaoWang2019}.  Gauss quadrature collocation methods use either Legendre-Gauss (LG) points \cite{Benson2,Garg1,Garg2,HagerHouRao15a,HagerHouRao16a,HagerLiuMohapatraWangRao18a,DuChenHager2019}, Legendre-Gauss-Radau (LGR) points \cite{Kameswaran1,Garg1,Garg2,Garg3,HagerHouRao15a,DuChenHager2019,HagerHouMohapatraRaoWang2019}, or Legendre-Gauss-Lobatto (LGL) points \cite{Elnagar1}. 

Various $h$ or $p$ direct collocation methods have been developed previously.  Reference~\cite{Gong3} describes what is essentially a $p$ method where a differentiation matrix is used to identify switches, kinks, corners, and other discontinuities in the solution.  Reference~\cite{Zhao2} develops a fixed-order method that uses a density function to generate a sequence of non-decreasing size meshes on which to solve the optimal control problem.   Finally, in Ref.~\cite{Betts3} an error estimate is developed by integrating the difference between an interpolation of the time derivative of the state and the right-hand side of the dynamics.  The error estimate developed in Ref.~\cite{Betts3} is predicated on the use of a fixed-order method (for example, trapezoid, Hermite-Simpson, Runge-Kutta) and computes a low-order approximation of the integral of the aforementioned difference.

Although $h$ methods have been used extensively and $p$ methods are useful on certain types of problems, both the $h$ and $p$ approaches have limitations.  In the case of an $h$ method, it may be required to use an extremely fine mesh to improve accuracy.  In the case of a $p$ method, it may be required to use an unreasonably large degree polynomial to improve accuracy.  In order to reduce significantly the size of the finite-dimensional approximation, and thus improve computational efficiency of solving the NLP, in recent years the new class of $hp$ collocation methods has been developed for solving an optimal control problem.  In an $hp$ method, both the number of mesh intervals and the degree of the approximating polynomial within each mesh interval are allowed to vary.  While $hp$ methods were originally developed as finite-element methods for solving partial differential equations \cite{Babuska1,Babuska2,Babuska3,Babuska4,Babuska5}, over the past several years $hp$ methods have been developed for solving optimal control problems \cite{Darby2,Darby3,Patterson2015,Liu2015,Liu2018}.  The methods described in Refs.~\cite{Darby2} and \cite{Darby3} describe methods where the error estimate is based on the difference between an approximation of the time derivative of the state and the right-hand side of the dynamics midway between the collocation points.  Next, Ref.~\cite{Patterson2015} develops an error estimate based on the difference between the state interpolated using an increased number of Legendre-Gauss-Radau points in each mesh interval and the state obtained by integrating the dynamics on the solution using the interpolated state and control.  Similar to the methods of Refs.~\cite{Darby2} and \cite{Darby3}, however, the method of Ref.~\cite{Patterson2015} can only increase the size of the mesh.  Additionally, Ref.~\cite{Liu2015} develops a method that adjusts the number of mesh intervals and the degree of the approximating polynomial within a mesh interval based on a proven convergence rate of the $hp$-adaptive LGR collocation method.  Moreover, Ref.~\cite{Liu2015} describes an approach for reducing the size of the mesh.  Finally, Ref.~\cite{Liu2018} describes an approach for adjusting the mesh based on the decay rate of the coefficients of a Legendre polynomial expansion of the state.

While the aforementioned $hp$ mesh refinement methods as described in Refs.~\cite{Darby2,Darby3,Patterson2015,Liu2015,Liu2018} can improve accuracy and computational efficiency when compared with traditional $h$ or $p$ methods, a key missing aspect of these methods is that they do not exploit the structure in the optimal solution.  In particular, in the case where the optimal solution is nonsmooth, rather than refining the mesh based on knowledge of the structure of the solution, these previously developed methods improve accuracy by increasing the number of collocation points in the vicinity of a discontinuity.  As a result, these methods often place an undesirably large number of collocation points in the neighborhood of a discontinuity.  Consequently, optimal control problems whose controls are discontinuous can require a large number of mesh refinement iterations to converge using the previously developed $hp$-adaptive methods.

A particular class of optimal control problems whose optimal solutions are discontinuous is the class of {\em bang-bang optimal control problems}.  Bang-bang optimal control problems arise in a wide variety of well-known application areas \cite{Shamsi1}.  A key feature of an optimal control problem with a bang-bang optimal control is that the Hamiltonian is linear with respect to one or more components of the control.  Due to the linear dependence of the Hamiltonian on the control and under the assumption that the solution does not contain any singular arcs, Pontryagin's minimum principle applies such that the optimal control lies at either its minimum or maximum limit.  In the context of a mesh refinement method using collocation, if the method could algorithmically detect the bang-bang structure of the optimal solution, it may be possible to obtain an accurate solution in a more computationally efficient manner than would be possible using a standard mesh refinement method that does not exploit the structure of the solution.  

The objective of this research is to develop a mesh refinement method for solving bang-bang optimal control problems by algorithmically exploiting the structure of the optimal solution.  In particular, this research focuses on the development of a method that significantly improves computational efficiency while simultaneously reducing the mesh size and the number of mesh refinement iterations required in order to obtain a solution to a bang-bang optimal control problem.  Previous research for solving bang-bang optimal control problems has been conducted using indirect methods.  Refs.~\cite{Maurer2,Kim1,Meier1,Hu1,Bertrand1,Huang1,Ledzewicz1,Maurer1} employ indirect shooting while Refs.~\cite{Ledzewicz1,Maurer1} include the second-order optimality conditions for bang-bang optimal control problems.  Furthermore, Refs.~\cite{Kaya1,Kaya2,Kaya3} employ a direct shooting method for solving bang-bang optimal control problems by parameterizing the control using piecewise constants where the durations of the intervals are added as optimization parameters in order to solve for the bang-bang control profile.  In particular, Ref.~\cite{Kaya3} combines the switch time computation and the time-optimal switching method developed in Refs.~\cite{Kaya1} and \cite{Kaya2}, respectively.  In the context of $hp$-adaptive mesh refinement methods, knotting methods have been developed in Refs.~\cite{Cuthrell1,Ross1} which allow discontinuities in optimal control profiles to be taken into account by explicitly introducing a switch time variable to the problem definition.  In particular, Ref.~\cite{Cuthrell1} utilizes the concept of super-elements which introduce switch time variables that are taken into account as parameters in the resulting optimization process.  The modified Legendre pseudospectral scheme described in Ref.~\cite{Shamsi1} also uses a similar approach by handling bang-bang optimal control problems using a knotting method that solves for the optimal control using an assumed number of switch times and constant control arcs.  The assumed number of switch times is then increased or decreased for each iteration based on the approximated solution, with the scheme converging upon the solution for the optimal number of switch times.  Additionally, Ref.~\cite{Agamawi2017} describes a slight modification of the $hp$-adaptive mesh refinement derived in Ref.~\cite{Liu2018}, where new mesh points are placed at discontinuity locations that are estimated based on the switching functions of the Hamiltonian of the optimal control problem.  Finally, Ref.~\cite{Miller2018} describes a mesh refinement method that is used in conjunction with the aforementioned $hp$-adaptive mesh refinement methods and detects discontinuities via jump function approximations.

In this paper a new direct collocation mesh refinement method is developed for solving optimal control problems whose solutions have a bang-bang structure.  While in principle the method of this paper can be used with any collocation method, in this paper the Legendre-Gauss-Radau collocation method \cite{Kameswaran1,Garg1,Garg2,Garg3,HagerHouRao15a,DuChenHager2019,HagerHouMohapatraRaoWang2019} is employed because it produces accurate state, control, and costate approximations \cite{Kameswaran1,Garg1,Garg2,Garg3,HagerHouMohapatraRaoWang2019}.  Moreover, the approach of this paper is designed to be used in conjunction with a previously developed $hp$ mesh refinement method such as those described in Refs.~\cite{Darby2,Darby3,Patterson2015,Liu2015,Liu2018}.   First, a solution is obtained on a coarse mesh.  Next, using the solution on this coarse mesh, the costate is estimated at the collocation points using the methods developed in Refs.~\cite{Garg1,Garg2,Garg3}.  Then, the state and costate approximations on the coarse mesh are used to determine algorithmically if the Hamiltonian is linear with respect to one or more components of the control.  Using the state and costate approximations, the switching functions are estimated at the collocation points for those components of the control for which the Hamiltonian depends upon linearly.  The estimates of the switching functions are then used to estimate a discontinuity in the control between any two collocation points where a switching function changes sign.  Using these estimates of the control discontinuities, the locations of the switch times in the control are then introduced as optimization variables and the optimal control problem is divided into multiple domains.  Within each domain, those components of the control that have a bang-bang solution structure are fixed at either their lower or upper limits depending upon the sign of the switching function.  It is important to note that those control components that do not have a bang-bang structure remain free to vary within their defined bounds.  Furthermore, it is noted that the multiple-domain formulation is analogous to using super-elements as described in Ref.~\cite{Cuthrell1}, with a key difference being that the structure of the bang-bang optimal control profile has been algorithmically detected using the estimates of the switching functions found on the initial mesh.  The multiple-domain optimal control problem is then solved using LGR collocation where it is noted again that the switch times are determined as part of the optimization.

This paper is organized as follows.  Section~\ref{sect:single-phase} defines the general single-phase optimal control problem in Bolza form.  Section~\ref{sect:LGR} describes the rationale for using Legend-Gauss-Radau collocation points as the set of nodes to discretize the continuous optimal control problem and the Legendre-Gauss-Radau collocation method.  Section~\ref{sect:Linear-Hamiltonian} briefly overviews the form of bang-bang optimal control problems.  Section~\ref{sect:bang-bang-method} describes the bang-bang mesh refinement method of this paper.  Section~\ref{sect:Examples} demonstrates the performance of the mesh refinement method of this paper when solving bang-bang optimal control problems as compared to the four previously developed $hp$-adaptive mesh refinement methods of Refs.~\cite{Patterson2015,Darby2,Liu2015,Liu2018}.  Section~\ref{sect:Discussion} provides a discussion of both the approach and the results.  Finally, Section~\ref{sect:Conclusions} provides conclusions on this research.

\section{Single-Phase Optimal Control Problem\label{sect:single-phase}} 

Without loss of generality, consider the following general single-phase optimal control problem in Bolza form defined on the time horizon $t\in[t_0,t_f]$.  Determine the state $\m{y}(t)\in\bb{R}^{1~\times~n_y}$, the control $\m{u}(t)\in\bb{R}^{1~\times~n_u}$,  the start time $t_0\in\bb{R}$, and the terminus time $t_f\in\bb{R}$ that minimize the objective functional
\begin{equation}\label{eq:single-cost-t}
  \C{J}=\C{M}(\m{y}(t_0),t_0,\m{y}(t_f),t_f)+ \int_{t_0}^{t_f} \C{L}(\m{y}(t),\m{u}(t), t)~dt~,
\end{equation}
subject to the dynamic constraints
\begin{equation}\label{eq:single-dyn-t}
  \frac{d\m{y}}{dt} = \m{a}(\m{y}(t),\m{u}(t), t)~,
\end{equation}
the inequality path constraints
\begin{equation}\label{eq:single-path-t}
\m{c}_{\min} \leq \m{c}(\m{y}(t),\m{u}(t), t)\leq \m{c}_{\max}~,
\end{equation}
and the boundary conditions
\begin{equation}\label{eq:single-bc-t}
  \m{b}_{\min} \leq \m{b}(\m{y}(t_0),t_0,\m{y}(t_f),t_f) \leq \m{b}_{\max}~.
\end{equation}

\section{Legendre-Gauss-Radau Collocation \label{sect:LGR}}

In order to develop the method described in this paper, a direct collocation method must be chosen.  While in principle any direct collocation method can be used to approximate the optimal control problem given in Section \ref{sect:single-phase}, in this research the previously developed Legendre-Gauss-Radau (LGR) collocation method \cite{Kameswaran1,Garg1,Garg2,Garg3,Patterson2015,HagerHouMohapatraRaoWang2019} will be used because it has been shown that the LGR collocation method produces an accurate state, control, and costate \cite{Kameswaran1,Garg1,Garg2,Garg3,HagerHouMohapatraRaoWang2019}.  It is noted that the accuracy of the costate estimate of the LGR collocation method plays an important role in the method of this paper because the accuracy of the costate directly influences the accuracy of the estimates of the switching functions of the Hamiltonian which are in turn used to estimate the switch times in the optimal control (forming the basis of the method of this paper).

In the context of this research, a multiple-interval form of the LGR collocation method is chosen.  The time horizon $t\in[t_0,t_f]$ may be divided into $Q$ time domains, $\C{P}_d=[t_s^{[d-1]},t_s^{[d]}]\subseteq[t_0,t_f],~d\in\{1,\ldots,Q\}$, such that
\begin{equation}\label{eq:time-domain-properties}
\bigcup_{d=1}^{Q}\C{P}_d=[t_0,t_f] ~, \quad\bigcap_{d=1}^{Q}\C{P}_d=\{t_s^{[1]},\ldots,t_s^{[Q-1]}\}~,
\end{equation}
where $t_s^{[d]},~d\in\{1,\ldots,Q-1\}$ are the switch time variables of the problem, $t_s^{[0]}=t_0$, and $t_s^{[Q]} = t_f$.  Thus in the case where $Q=1$ the domain consists of only a single domain $\C{P}_1=[t_0,t_f]$ and $\{t_s^{[1]},\ldots,t_s^{[Q-1]}\}=\emptyset$.  In the multiple-interval LGR collocation method, each of the time domains $\C{P}_d=[t_s^{[d-1]},t_s^{[d]}],~d\in\{1,\ldots,Q\}$, is converted into the domain $\tau\in[-1,+1]$ using the affine transformation,
\begin{equation}\label{eq:affine-transformation}
  \begin{array}{lcl}
    t & = & \displaystyle \frac{t_s^{[d]}-t_s^{[d-1]}}{2}\tau + \frac{t_s^{[d]}+t_s^{[d-1]}}{2}~, \vspace{3pt} \\
    \tau & = &  \displaystyle 2\frac{t-t_s^{[d-1]}}{t_s^{[d]}-t_s^{[d-1]}}-1~.
  \end{array}
\end{equation}
The interval $\tau\in[-1,+1]$ for each domain $\C{P}_d$ is then divided into $K$ mesh intervals, $\C{S}_k=[T_{k-1},T_k]\subseteq [-1,+1],\;k\in\{1,\ldots,K\}$ such that
\begin{equation}\label{eq:mesh-interval-properties}
\bigcup_{k=1}^K\C{S}_k=[-1,+1] ~, \quad\bigcap_{k=1}^K\C{S}_k=\{T_1,\ldots,T_{K-1}\}~,
\end{equation}
and $-1=T_0<T_1<\ldots<T_{K-1}<T_K=+1$.  For each mesh interval, the LGR points used for collocation are defined in the domain of $[T_{k-1},T_k]$ for $k\in\{1,\ldots,K\}$.  The state of the continuous optimal control problem is then approximated in mesh interval $\C{S}_k,\;k\in\{1,\ldots,K\}$, as 
\begin{equation}\label{eq:LGR-state-approximation}
\m{y}^{(k)}(\tau)  \approx \m{Y}^{(k)}(\tau) = \sum_{j=1}^{N_k+1} \m{Y}_{j}^{(k)}
\ell_{j}^{(k)}(\tau)~, \quad  \ell_{j}^{(k)}(\tau) = \prod_{\stackrel{l=1}{l\neq j}}^{N_k+1}\frac{\tau-\tau_{l}^{(k)}}{\tau_{j}^{(k)}-\tau_{l}^{(k)}}~, 
\end{equation}  
where
$\ell_{j}^{(k)}(\tau)$ for $ j\in\{1,\ldots,N_k+1\}$ is a basis of Lagrange polynomials on $\C{S}_k$, $\left(\tau_1^{(k)},\ldots,\tau_{N_k}^{(k)}\right)$ are the set of $N_k$ Legendre-Gauss-Radau (LGR) \cite{Abramowitz1} collocation points in the interval $[T_{k-1},T_k)$, $\tau_{N_k+1}^{(k)}=T_k$ is a non-collocated support point, and $\m{Y}_{j}^{(k)} \equiv  \m{Y}^{(k)}(\tau_j^{(k)})$.  Differentiating $\m{Y}^{(k)}(\tau)$ in Eq.~(\ref{eq:LGR-state-approximation}) with respect to $\tau$ gives
\begin{equation}\label{eq:LGR-diff-state-approximation}
  \frac{d\m{Y}^{(k)}(\tau)}{d\tau} = \sum_{j=1}^{N_k+1}\m{Y}_{j}^{(k)}\frac{d\ell_j^{(k)}(\tau)}{d\tau}~.
\end{equation}
The dynamics are then approximated at the $N_k$ LGR points in mesh
interval $k\in\{1,\ldots,K\}$ as
\begin{equation}\label{eq:LGR-defect}
 \sum_{j=1}^{N_k+1}D_{lj}^{(k)} \m{Y}_j^{(k)} - \frac{t_f-t_0}{2}\m{a}\left(\m{Y}_l^{(k)},\m{U}_l^{(k)},t (\tau_l^{(k)},t_0,t_f)\right) = \m{0} ~,\quad l \in \{1,\ldots,N_k\}~,
\end{equation}
where 
\begin{equation*}
  D_{lj}^{(k)} = \frac{d\ell_j^{(k)}(\tau_l^{(k)})}{d\tau}~,\quad l \in \{1,\ldots,N_k\}~,~ j \in \{1,\ldots,N_k+1\}~,
\end{equation*}
are the elements of the $N_k\times (N_k+1)$ {\em Legendre-Gauss-Radau differentiation matrix} \cite{Garg1} in mesh interval $\C{S}_k$, $\;k\in\{1,\ldots,K\}$, and $\m{U}_l^{(k)}$ is the approximation of the control at the $l^{th}$ collocation point in mesh interval $\C{S}_k$.  It is noted that continuity in the state and time between mesh intervals $\C{S}_{k-1}$ and $\C{S}_{k}$, $k\in\{1,\ldots,K\}$, is enforced by using the same variables to represent $\m{Y}_{N_{k-1}+1}^{(k-1)}$ and $\m{Y}_{1}^{(k)}$, while continuity in the state between the domains $\C{P}_{d-1}$ and $\C{P}_{d}$, $d\in\{2,\ldots,Q\}$, is achieved using the additional continuity constraint
\begin{equation}\label{eq:domain-continuity}
\m{Y}_{N^{[d-1]}+1}^{[d-1]} = \m{Y}_{1}^{[d]}~,
\end{equation}
where the superscript $[d]$ is used to denote the $d^{th}$ time domain, $\m{Y}_{j}^{[d]}$ denotes the value of the state approximation at the $j^{th}$ discretization point in the time domain $\C{P}_d$, and $N^{[d]}$ is the total number of collocation points used in time domain $\C{P}_d$ computed by
\begin{equation}\label{eq:N-sum}
N^{[d]} = \sum_{k=1}^{K^{[d]}} N_k^{[d]}~.
\end{equation}

The Legendre-Gauss-Radau approximation of the multiple-domain optimal control problem then leads to the following nonlinear programming problem (NLP).  Minimize the objective function
\begin{equation}\label{eq:NLP-cost}
  \C{J}=\C{M}(\m{Y}_1^{[1]},t_0,\m{Y}_{N^{[Q]}+1}^{[Q]},t_f)+ \sum_{d=1}^{Q} \frac{t_s^{[d]}-t_s^{[d-1]}}{2}\left[\m{w}^{[d]}\right]\tr\m{L}^{[d]}~,
\end{equation}
subject to the defect constraints
\begin{equation}\label{eq:NLP-defect}
\g{\Delta}^{[d]} =
\m{D}^{[d]}\m{Y}^{[d]} - \frac{t_s^{[d]}-t_s^{[d-1]}}{2}\m{A}^{[d]}
=\m{0}
~,\quad  d \in \{1,\ldots,Q\}~,
\end{equation}
the path constraints
\begin{equation}\label{eq:NLP-path}
\m{c}_{\min} \leq
\m{C}_{j}^{[d]}
\leq \m{c}_{\max}
~,\quad  j \in \{1,\ldots,N^{[d]}\}~,~ d \in \{1,\ldots,Q\}~,
\end{equation}
the boundary conditions
\begin{equation}\label{eq:NLP-bc}
  \m{b}_{\min} \leq \m{b}(\m{Y}_1^{[1]},t_0,\m{Y}_{N^{[Q]}+1}^{[Q]},t_f) \leq \m{b}_{\max}~.
\end{equation}
and the continuity constraints
\begin{equation}\label{eq:NLP-continuity}
\m{Y}_{N^{[d-1]}+1}^{[d-1]} = \m{Y}_{1}^{[d]}~,~ d \in \{2,\ldots,Q\}~,
\end{equation}
where
\begin{equation}\label{eq:NLP-A-def}
\m{A}^{[d]} = \begin{bmatrix}
\m{a}\left(\m{Y}_{1}^{[d]},\m{U}_{1}^{[d]},t_1^{[d]}\right)\\
\vdots\\
\m{a}\left(\m{Y}_{N^{[d]}}^{[d]},\m{U}_{N^{[d]}}^{[d]},t_{N^{[d]}}^{[d]}\right)
\end{bmatrix}
\in \mathbb{R}^{N^{[d]}~\times~n_y}~,
\end{equation}
\begin{equation}\label{eq:NLP-C-def}
\m{C}^{[d]} = \begin{bmatrix}
\m{c}\left(\m{Y}_{1}^{[d]},\m{U}_{1}^{[d]},t_1^{[d]}\right)\\
\vdots\\
\m{c}\left(\m{Y}_{N^{[d]}}^{[d]},\m{U}_{N^{[d]}}^{[d]},t_{N^{[d]}}^{[d]}\right)
\end{bmatrix}
\in \mathbb{R}^{N^{[d]}~\times~n_c}~,
\end{equation}
\begin{equation}\label{eq:NLP-G-def}
\m{L}^{[d]} = \begin{bmatrix}
\C{L}\left(\m{Y}_{1}^{[d]},\m{U}_{1}^{[d]},t_1^{[d]}\right)\\
\vdots\\
\C{L}\left(\m{Y}_{N^{[d]}}^{[d]},\m{U}_{N^{[d]}}^{[d]},t_{N^{[d]}}^{[d]}\right)
\end{bmatrix}
\in \mathbb{R}^{N^{[d]}~\times~1}~,
\end{equation}
$\m{D}^{[d]} \in \mathbb{R}^{N^{[d]}~\times~[N^{[d]}+1]}$ is the LGR differentiation matrix in time domain $\C{P}_d,~d \in\{1,\ldots,Q\}$, and $\m{w}^{[d]} \in \mathbb{R}^{N^{[d]}~\times~1}$ are the LGR weights at each node in time domain $\C{P}_d,~ d \in\{1,\ldots,Q\}$.  It is noted that $\m{a} \in \mathbb{R}^{1~\times~n_y}$, $\m{c} \in \mathbb{R}^{1~\times~n_c}$, and $\C{L} \in \mathbb{R}^{1~\times~1}$ correspond, respectively, to the vector fields that define the right-hand side of the dynamics, the path constraints, and the integrand of the optimal control problem, where $n_y$ and $n_c$ are, respectively, the number of state components and path constraints in the problem.  Additionally, the state matrix, $\m{Y}^{[d]}\in \mathbb{R}^{[N^{[d]}+1]~\times~n_y}$, and the control matrix, $\m{U}^{[d]}\in \mathbb{R}^{N^{[d]}~\times~n_u}$, in time domain $\C{P}_d,~ d \in \{1,\ldots,Q\}$, are formed as 
\begin{equation}\label{eq:NLP-Var-Mat}
\m{Y}^{[d]} = \begin{bmatrix}
\m{Y}^{[d]}_{1}\\
\vdots \\
\m{Y}^{[d]}_{N^{[d]}+1}
\end{bmatrix}
\text{ and }
\m{U}^{[d]} = 
\begin{bmatrix}
\m{U}^{[d]}_{1}  \\
\vdots	\\
\m{U}^{[d]}_{N^{[d]}}
\end{bmatrix}~,
\end{equation}
respectively, where $n_u$ is the number of control components in the problem.  Finally, as described in Ref.~\cite{Garg1}, estimates of the costate may be obtained at each of the discretization points in the time domain $\C{P}_d, d\in\{1,\ldots,Q\}$, using the transformation
\begin{equation}\label{eq:NLP-costates}
\g{\lambda}^{[d]} = (\m{W}^{[d]})^{-1}\g{\Lambda}^{[d]}~,\quad \g{\lambda}_{N^{[d]}+1}^{[d]} = (\m{D}_{N^{[d]}+1}^{[d]})\tr \g{\Lambda}^{[d]}~,
\end{equation}
where $\g{\lambda}^{[d]}\in\mathbb{R}^{N^{[d]}~\times~n_y}$ is a matrix of the costate estimates at the collocation points in time domain $\C{P}_d$, $\m{W}^{[d]} = \text{diag}(\m{w}^{[d]})$ is a diagonal matrix of the LGR weights at the collocation points in time domain $\C{P}_d$, $\g{\Lambda}^{[d]}\in\mathbb{R}^{N^{[d]}~\times~n_y}$ is a matrix of the NLP multipliers obtained from the NLP solver corresponding to the defect constraints at the collocation points in time domain $\C{P}_d$ , $\g{\lambda}_{N^{[d]}+1}^{[d]}\in\mathbb{R}^{1~\times~n_y}$ is a row vector of the costate estimates at the non-collocated end point in time domain $\C{P}_d$, and $\m{D}_{N^{[d]}+1}^{[d]}\in\mathbb{R}^{N^{[d]}~\times~1}$ is the last column of the LGR differentiation matrix in time domain $\C{P}_d$.


\section{Control-Linear Hamiltonian  \label{sect:Linear-Hamiltonian}}

The Hamiltonian of the Bolza optimal control problem defined in Eqs.~\eqref{eq:single-cost-t}--\eqref{eq:single-bc-t} is given as
\begin{equation}\label{Hamiltonian}
  \C{H}(\m{y}(t),\g{\lambda}(t),\m{u}(t),t) =\C{L}(\m{y}(t),\m{u}(t), t) + \g{\lambda}(t)\m{a}\tr(\m{y}(t),\m{u}(t),t).  
\end{equation}
Assume now that the Hamiltonian given in Eq.~\eqref{Hamiltonian} has the following form:
\begin{equation}\label{control-linear-Hamiltonian}
  \C{H}(\m{y}(t),\g{\lambda}(t),\m{u}(t),t) = f(\m{y}(t),\g{\lambda}(t),\m{z}(t),t) + \g{\sigma}(\m{y}(t),\g{\lambda}(t),t)\m{v}\tr(t),  
\end{equation} 
where $f\in\mathbb{R}$, $\g{\sigma}\in\mathbb{R}^{1~\times~I}$, $\m{v}(t)=(u_1(t),\ldots,u_I(t))\in\mathbb{R}^{1\times I}$ and $\m{z}(t)=(u_{I+1}(t),\ldots,u_{n_u}(t))\in\mathbb{R}^{1\times (n_u-I)}$.  It is noted that $\m{v}(t)$ and $\m{z}(t)$ are vectors that correspond, respectively, to those components of the control upon which the Hamiltonian depends linearly and nonlinearly, while $\g{\sigma}=(\sigma_1,\ldots,\sigma_I)$ is a vector that defines the $I$ switching functions.  In other words, the components $(u_1(t),\ldots,u_{n_u}(t))$ of the control are ordered such that any control-linear component has an index that is lower than any control-nonlinear component and $(\sigma_1,\ldots,\sigma_I)$ are functions of the state and the costate that are used to determine locations where a switch in the control may occur.  Then, assuming that the optimal solution contains no singular arcs, the optimal control for any control-linear component, $u_i^*(t),~i\in\{1,\ldots,I\}$, is obtained from Pontryagin's minimum principle as
\begin{equation}\label{eq:Pontryagin} 
  u_i^*(t) = \left\{
    \begin{array}{lcl}
      u_i^{\min} & , & \sigma_i(\m{y}(t),\g{\lambda}(t),t)>0, \\
      u_i^{\max} & , & \sigma_i(\m{y}(t),\g{\lambda}(t),t)<0,
    \end{array}
    \right.\quad (i=1,\ldots,I)~,
\end{equation}
where the signs of the switching functions, $\sigma_i,~i\in\{1,\ldots,I\}$, are determined by the state $\m{y}(t)$, the costate $\g{\lambda}(t)$, and the time $t$.  Thus, the optimal control will have a bang-bang structure where discontinuities in the solution occur whenever a switching function, $\sigma_i,~i\in\{1,\ldots,I\}$, changes sign (again, assuming no singular arcs).

Assume now that the optimal control problem has been approximated using the LGR collocation method as described in Section \ref{sect:LGR}.  Furthermore, assume that the NLP resulting from LGR collocation has been solved to obtain estimates of the state, control, and costate as given in Eqs.~\eqref{eq:NLP-Var-Mat} and \eqref{eq:NLP-costates}.  Then the approximation can be utilized to detect possible switch times in the optimal control.  These possible discontinuity locations can be detected by evaluating the switching functions $\sigma_i,~i\in\{1,\ldots,I\}$, at each of the collocation points and checking for sign changes in $\sigma_i$ between two adjacent collocation points.  In particular, a sign change in $\sigma_i$ between any two adjacent collocation points $t_k$ and $t_{k+1}$ indicates that $\sigma_i$ has a root in the time interval $t\in[t_k,t_{k+1}]$, thus indicating a switch in the control at a time $t\in[t_k,t_{k+1}]$.  As discussed in Ref.~\cite{bang_bang}, the rationale for using changes in the sign of $\sigma_i$ to detect switches in the control is due to the fact that the gradient of the objective with respect to a switching point in a bang-bang control is proportional to the coefficient $\sigma_i$ of the control in the Hamiltonian.  Consequently, the gradient of the objective with respect to a switching point only vanishes at values of time where the coefficient of the control vanishes, and the change in sign of $\sigma_i$ between a pair of points provides an estimate of a switch time in the component of the control corresponding to the switching function $\sigma_i$.

\section{Bang-Bang Mesh Control Refinement Method \label{sect:bang-bang-method}}

In this section, the bang-bang mesh refinement method is developed.  The method consists of several steps.  The first step, described in Section \ref{sect:bang-bang-detection} determines if the Hamiltonian is linear in one or more components of the control.  The second step, described in Section \ref{sect:structure-detection} computes estimates of the locations of discontinuities in each component of the control upon which the Hamiltonian depends upon linearly.  The third step, described in Section \ref{sect:variable-mesh-points}, reformulates the optimal control problem into a multiple-domain optimal control problem where the domains are determined based on the estimates of the discontinuity locations obtained in the second step of the method.  The interior endpoints of these domains are then treated as optimization variables by introducing an appropriate number of switch time parameters which use the estimates of the discontinuity locations as initial guesses for these new variables.  In addition, any component of the control for which the Hamiltonian has been determined to depend upon linearly is then fixed at either its minimum or maximum value in each domain depending upon the direction of the switch in sign of the corresponding switching function component identified in the second step.  The multiple-domain optimal control problem is then solved iteratively using the LGR collocation method described in Section \ref{sect:LGR} together with a previously developed $hp$-adaptive mesh refinement method (to refine the mesh based on the error in the solution for those control components for which the Hamiltonian does not depend upon linearly and provide sufficient resolution in the state approximation).  Finally, Section \ref{sect:summary} provides a summary of the bang-bang mesh refinement method.  

\subsection{Method for Identifying Bang-Bang Optimal Control Problems \label{sect:bang-bang-detection}}

Assume now that the optimal control problem as formulated in Section~\ref{sect:single-phase} has been transcribed into a nonlinear programming problem using the LGR collocation method as developed in Section~\ref{sect:LGR} assuming that $Q=1$ (that is, a single domain is used).  Furthermore, assume that the resulting NLP has been solved to obtain estimates of the state, control, and costate as given in Eqs.~\eqref{eq:NLP-Var-Mat} and \eqref{eq:NLP-costates}, and that the mesh refinement accuracy tolerance has not been satisfied and, thus, mesh refinement is required.  As a result, it is possible that the optimal solution may possess a bang-bang optimal control.  The first step in determining if a bang-bang optimal control is a possibility is to determine if the Hamiltonian is linear in one or more of the components of the control.  Moreover, the method should be able to detect a control-linear Hamiltonian without requiring any external intervention.

In this research, the determination of a control-linear Hamiltonian is made using hyper-dual derivative approximations~\cite{Fike2011}.  In particular, using the state and costate approximations obtained from the LGR collocation method, the hyper-dual derivative approximation of Ref.~\cite{Fike2011} is used to compute the first- and second-derivatives of the Hamiltonian with respect to each component of the control.  In order to identify linearity in the Hamiltonian with respect to a given control component, sample values of the control component are taken between the bounds of the control and partial derivatives of the Hamiltonian are then taken with respect to that control component while holding all other variables constant.  If the partial derivatives of the Hamiltonian obtained using the various sample values of the control component are found to be constant, then the Hamiltonian is identified as possibly being linear with respect to that control component.  Moreover, if any control component has been identified to have a zero second derivative and has zero second-order partial derivatives with respect to any other control components [that is, all cross partial derivatives $\frac{\partial^2 H}{\partial u_i\partial u_j}=0,~(i\neq j)$], then the Hamiltonian is linear with respect to that component of the control and is thus a candidate for bang-bang control mesh refinement.  Estimates of the switching functions of the Hamiltonian are then obtained by computing the partial of the Hamiltonian with respect to each control-linear component using hyper-dual derivative approximations.  Finally, it is noted that the hyper-dual derivative approximation is not subject to truncation error up to second derivatives and, thus, provides exact (that is, to machine precision) first- and second-derivatives of the Hamiltonian with respect to the control.  
 
\subsection{Estimating Locations of Switches in Control\label{sect:structure-detection}}

Once it has been determined that the Hamiltonian is linear in at least one control component, the next step is to estimate times at which a control component may switch between its lower and upper limit, thus leading to a discontinuity in the control solution profile.  Assuming that the optimal solution contains no singular arcs, a discontinuity in the control will occur when a switching function $\sigma_i,~i\in\{1,\ldots,I\}$ changes sign.  Given that the solution of the optimal control problem has been approximated using the LGR collocation as given in Section \ref{sect:LGR}, an {\em estimate} of a discontinuity in a control component will be when a switching function $\sigma_i$ changes sign between two adjacent collocation points.  Furthermore, because any switching function, $\sigma_i$ may change sign one or more times, it is possible that any or all components of the control $\m{v}(t)$ may have one or more discontinuities and that the discontinuity time for a particular component, $v_i(t)$, may differ from the discontinuity time of another control-linear component, $v_j(t),~(i\neq j)$.

\begin{figure}[!ht]
\centering
\includegraphics[keepaspectratio=true,width=2.5in]{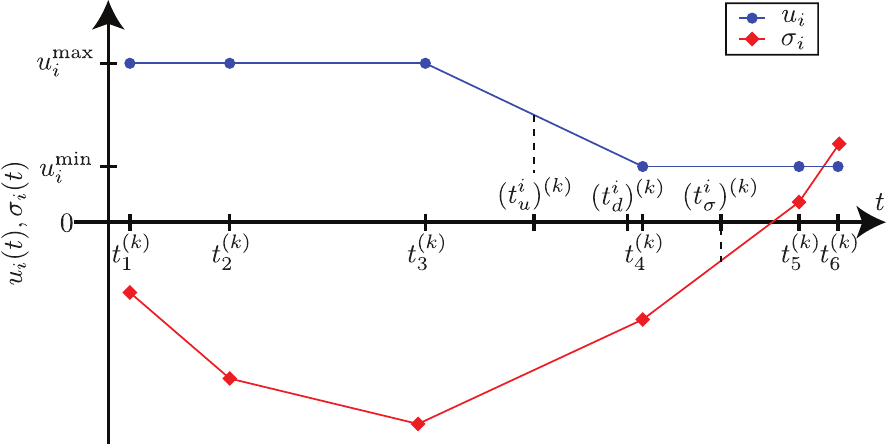}
\renewcommand{\baselinestretch}{1}\normalsize\normalfont
\caption{Estimates of discontinuity location $(t_d^{i})^{(k)}$ corresponding to control component $u_i$ in mesh interval $\C{S}_k$ using corresponding switching function $\sigma_i(t)$ with six collocation points $\{t_1^{(k)},\ldots,t_6^{(k)}\}$. \label{fig:DiscontinuityEstimation}}

\end{figure}

In order to accommodate the possibility of multiple discontinuities within a given mesh interval, $\C{S}_k,~k\in\{1,\ldots,K\}$, each discontinuity estimate time $(t_d^{i})^{(k)}$ in $\C{S}_k$ for a given control component $u_i,~i\in\{1,\ldots,I\}$ is computed using
\begin{equation}\label{eq:discontinuity-estimate-time}
(t_d^{i})^{(k)} = \frac{ (t_{\sigma}^i)^{(k)} + (t_{u}^i)^{(k)}}{2}~,
\end{equation}
where $(t_{\sigma}^i)^{(k)}$ is the midpoint time between the two adjacent discretization points showing a change in sign of the corresponding switching function $\sigma_i$ within $\C{S}_k$, and $(t_{u}^i)^{(k)}$ is the midpoint time between the two adjacent discretization points displaying the largest absolute difference in the value of the $i^{th}$ control-linear component $u_i$ within $\C{S}_k$.  An example schematic of estimating the discontinuity time, $(t_d^{i})^{(k)}$, for control-linear component, $u_i$, for mesh interval $\C{S}_k$ is shown in Fig.~\ref{fig:DiscontinuityEstimation}.

\begin{figure}[!ht]
\centering

\includegraphics[scale=1]{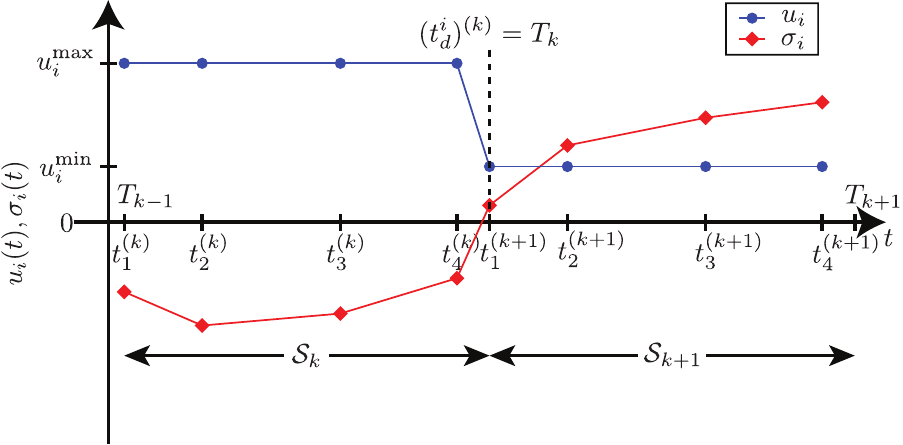}
\renewcommand{\baselinestretch}{1}\normalsize\normalfont
\caption{Estimates of discontinuity location $(t_d^{i})^{(k)}$ corresponding to control component $u_i$ across mesh intervals $\C{S}_k$ and $\C{S}_{k+1}$ using corresponding switching function $\sigma_i(t)$ and four collocation points in each mesh interval $\{t_1^{(k)},\ldots,t_4^{(k)},t_1^{(k+1)},\ldots,t_4^{(k+1)}\}$. \label{fig:DiscontinuityCaptured}}

\end{figure}

By taking into account both of the midpoint times $(t_{\sigma}^i)^{(k)}$ and $(t_{u}^i)^{(k)}$ in estimating the discontinuity time $(t_d^i)^{(k)}$, the bang-bang control profile can be properly maintained relative to all control-linear components within a mesh interval $\C{S}_k$ containing multiple discontinuities.  Furthermore, in the event that the switching function $\sigma_i$ changes sign across two adjacent mesh intervals $\C{S}_k$ and $\C{S}_{k+1}$, the mesh point $T_k$ that lies at the interface between the mesh intervals $\C{S}_k$ and $\C{S}_{k+1}$ is used as the estimate for the discontinuity time $(t_d^i)^{(k)}$ (that is $(t_d^i)^{(k)}=T_k$) as shown in Fig.~\ref{fig:DiscontinuityCaptured}.

After checking all mesh intervals $\C{S}_k,~k \in\{1,\ldots,K\}$, for possible discontinuities, the computed discontinuity estimates, $(t_d^i)^{(k)},~i\in\{1,\ldots,I\},~k\in\{1,\ldots,K\}$, are arranged in ascending order and used as initial guesses for the switch time parameters $t_s^{[S]},~S\in\{1,\ldots,n_s\}$, that are to be introduced in the subsequent mesh iterations, where $n_s$ is equal to the total number of discontinuities detected in the solution obtained on the initial mesh.  Finally, the limit of each control-linear component $u_i,~i\in\{1,\ldots,I\},$ to be used on the subsequent mesh iterations in each of the newly created time domains $\C{P}_d,~d\in\{1,\ldots,Q\}$, is identified by checking the sign of the corresponding switching function $\sigma_i$ within $\C{P}_d$, while each control-nonlinear component $u_i,~i\in\{I+1,\ldots,n_u\},$ is left free to vary between its defined bounds, where $Q = n_s + 1$ on the subsequent mesh iterations.  Thus, by using the estimates of the switching functions, the structure of the bang-bang optimal control profile has been automatically detected and used to set up the subsequent mesh iterations to include $n_s$ switch time parameters.

\subsection{Reformulation of Optimal Control Problem Into Multiple Domains\label{sect:variable-mesh-points}}

\begin{figure}[!ht]
\centering

\includegraphics[keepaspectratio=true,width=2.5in]{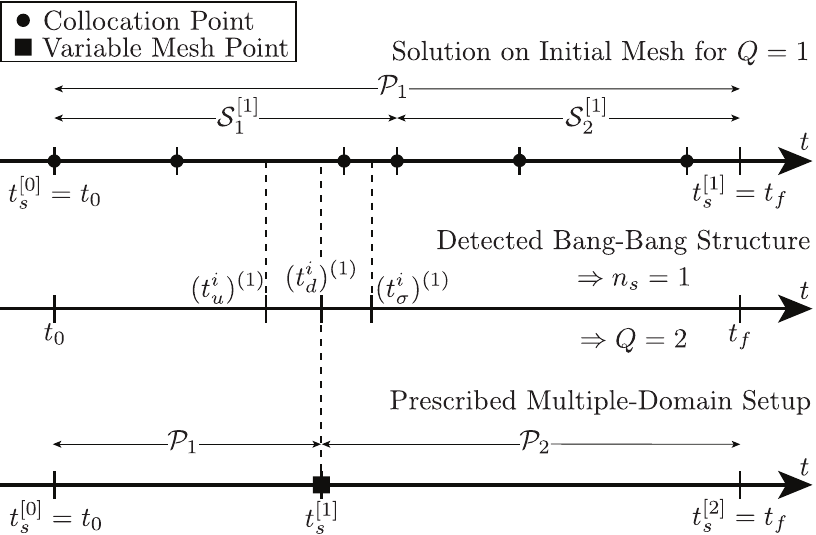}

\renewcommand{\baselinestretch}{1}\normalsize\normalfont
\caption{Schematic of process that creates a multiple-domain optimal control problem with $Q=n_s+1$ domains where variable mesh points are introduced as optimization parameters in order to determine the $n_s$ optimal switch times in the components of the control for which the Hamiltonian depends upon linearly. \label{fig:MeshRefinementSchematic}}

\end{figure}

Assuming the optimal control problem has been identified suitable for bang-bang mesh refinement as described in Section~\ref{sect:bang-bang-detection}, the method for automatically detecting the structure of the bang-bang control as described in Section~\ref{sect:structure-detection} may be employed.  Once acquired, the detected structure of the bang-bang control may be used to introduce the appropriate number of switch time parameters, $t_s^{[S]},~S \in\{1,\ldots,n_s\}$, to be solved for on subsequent mesh iterations, where the initial guess for each switch time parameter $t_s^{[S]}$ is the estimated discontinuity time that was found using the process described in Section~\ref{sect:structure-detection}.  The switch time parameters are included by using variable mesh points between the newly created time domains, $\C{P}_d,~d \in\{1,\ldots,Q\}$, where $Q=n_s+1$ .  Specifically, the variable mesh points are employed by dividing the time horizon $t=[t_0,t_f]$ of the optimal control problem identified as a bang-bang control into $Q=n_s+1$ time domains as described in Section~\ref{sect:LGR}, where $n_s$ is the number of switch time parameters introduced based on the detected structure of the bang-bang control profile for the solution on the initial mesh.  Each of the time domains $\C{P}_d=[t_s^{[d-1]},t_s^{[d]}],~d\in\{1,\ldots,Q\}$, has bounds enforced on $t_s^{[d-1]}$ and $t_s^{[d]}$ which appropriately bracket the estimated discontinuity times detected from the structure.  Additionally, the bounds of the control-linear components $u_i,~i \in\{1,\ldots,I\},$ within each time domain $\C{P}_d,~d \in\{1,\ldots,Q\}$, are set to the identified constant bang control limit for each control-linear component during the corresponding spans of time based on the previously detected structure, while the bounds of the control-nonlinear components  $u_i,~i \in\{I+1,\ldots,n_u\},$ are left unchanged.  Thus the bang-bang optimal control problem is effectively transcribed into a multiple-domain optimal control problem employing constant values for the control-linear components in each time domain $\C{P}_d~,d \in\{1,\ldots,Q\},$ such that the optimal switch times are solved for on the subsequent mesh iterations while the control-nonlinear components are left free to vary between their respective bounds.  A schematic for the process of dividing the identified bang-bang optimal control problem into a multiple-domain optimal control problem employing constant values for the control-linear components $u_i,~i \in\{1,\ldots,I\},$ in each time domain is shown in Fig.~\ref{fig:MeshRefinementSchematic}.

\subsection{Summary of Bang-Bang Control Mesh Refinement Method\label{sect:summary}}

A summary of the bang-bang control mesh refinement method developed in this paper appears below.  Here $M$ denotes the mesh number, and in each loop of the method, the mesh number increases by $1$.  The method terminates in Step~\ref{step:terminate} when the error tolerance is satisfied or when $M$ reaches a prescribed maximum $M_{\max}$.

\begin{shadedframe}
  \renewcommand{\baselinestretch}{1}\normalsize\normalfont
\vspace{-10pt}
\begin{center}
 \shadowbox{\bf Bang-Bang Control Mesh Refinement Method}
\end{center}
\begin{enumerate}[{\bf Step 1:}]
\item Specify initial mesh.
\item Solve LGR collocation NLP of Eqs.~\eqref{eq:NLP-cost}--\eqref{eq:NLP-continuity} on current mesh.\label{step:solve}
\item Compute relative error $e$ on current mesh.
\item If $e<\epsilon$ or $M>M_{\max}$, then quit. Otherwise, go to {\bf Step \ref{step:start-process}}. \label{step:terminate}
\item If $M=1$, determine the number $I$ of control components for which Hamiltonian is linear in the control in manner described in Section~\ref{sect:bang-bang-detection}. Otherwise, go to {\bf Step \ref{step:standard}}. \label{step:start-process}
\item If $I=0$ or $M>1$, employ standard mesh refinement and return to {\bf Step \ref{step:solve}}.\label{step:standard}
\item If $I>0$ and $M=1$, employ bang-bang mesh refinement using the following steps:
  \begin{enumerate}
  \item Estimate $n_s$ discontinuities in control components using method described in Section~\ref{sect:structure-detection}.
  \item Partition time domain into multiple ($Q=n_s+1$) domains using method of Section \ref{sect:variable-mesh-points}.
  \item Solve multiple-domain optimal control problem that includes following features:
    \begin{enumerate}[(i)]
      \item Include $n_s$ variable mesh points in multiple-domain formulation.
      \item Fix $I$ bang-bang control components at either lower/upper limit in each of the $Q$ domains.
      \end{enumerate}
    \item Increment $M\longrightarrow M+1$ and return to {\bf Step \ref{step:solve}}.  
  \end{enumerate}
\end{enumerate}
\end{shadedframe}

\section{Examples \label{sect:Examples}}

In this section, three nontrivial bang-bang optimal control problems are solved using the bang-bang mesh refinement method described in Section~\ref{sect:bang-bang-method}.  The first example is the three compartment model problem taken from Ref.~\cite{Ledzewicz1}.  The second example is the robot arm problem taken from Ref.~\cite{Munson}.  The third example is the free-flying robot problem taken from Ref.~\cite{Sakawa1}.  The efficiency of the bang-bang mesh refinement method developed in this paper is evaluated and compared against four previously developed $hp$-adaptive mesh refinement methods described in Refs.~\cite{Patterson2015,Darby2,Liu2015,Liu2018}.  For problems requiring more than a single mesh refinement to meet accuracy tolerance, the bang-bang mesh refinement method will utilize the $hp$-adaptive method described in Ref.~\cite{Liu2018} to further refine the phases not meeting the specified mesh accuracy tolerance.  It is noted that any of the four previously developed $hp$-adaptive mesh refinement methods may be paired with the bang-bang mesh refinement method.  An initial coarse mesh of ten intervals with five collocation points in each interval is used for each problem.  Furthermore, upon identification of the bang-bang control solution profile, the bang-bang mesh refinement method initially employs two mesh intervals to discretize each of the newly created time domains employing constant values for the control-linear components, with five collocation points in each mesh interval for the first and second examples and six collocation points in each mesh interval for the third example.  The performance of the mesh refinement methods are evaluated based on the number of mesh iterations, total number of collocation points used, and the total computation time for the problem to be solved satisfactorily for the specified mesh accuracy tolerance.  All plots were created using MATLAB Version R2016a (build 9.0.0.341360).

The following conventions are used for all of the examples.  First, $M$ denotes the number of mesh refinement iterations required to meet the mesh refinement accuracy tolerance of $\epsilon=10^{-6}$, where $M=1$ corresponds to the initial mesh.  Second, $N_f$ denotes the total number of collocation points on the final mesh.  Third, $\C{T}$ denotes the total computation time required to solve the optimal control problem.  Furthermore, the $hp$-adaptive mesh refinement methods described in Refs.~\cite{Patterson2015,Darby2,Liu2015,Liu2018} are referred to, respectively, as the $hp$-I, $hp$-II, $hp$-III, and $hp$-IV mesh refinement methods.  Additionally, the bang-bang mesh refinement method developed in this paper is referred to as the $hp$-BB mesh refinement method.  Finally, for all mesh refinement methods, a minimum and maximum number of  three and ten collocation points in each mesh interval is enforced, respectively.  All results shown in this paper were obtained using the C++ optimal control software $\mathbb{CGPOPS}$ \cite{AgamawiRaoCGPOPS2019} using the NLP solver IPOPT \cite{Biegler2} in full Newton (second-derivative) mode with an NLP solver optimality tolerance of $10^{-9}$.  All first- and second-derivatives required by the NLP solver were computed using the hyper-dual derivative approximation method described in Ref.~\cite{Fike2011}.  All computations were performed on a 2.9 GHz Intel Core i7 MacBook Pro running MAC OS-X version 10.13.6 (High Sierra) with 16GB 2133MHz LPDDR3 of RAM.  C++ files were compiled using Apple LLVM version 9.1.0 (clang-1000.10.44.2).

\subsection{Example 1: Three Compartment Model Problem \label{sect:threeCompartmentModel}}

Consider the following optimal control problem taken from Ref.~\cite{Ledzewicz1}.  Minimize the objective functional  
\begin{equation}\label{eq:threeCompartmentModel-Cost}
  \C{J} = r_1N_1(t)+r_2N_2(t)+r_3N_3(t)+\int_0^{T}u_1(t)~dt~,
\end{equation}
subject to the dynamic constraints
\begin{equation}\label{eq:threeCompartmentModel-Dynamics}
\begin{array}{l}
 \begin{array}{lcl}
   \dot{N}_1(t) & = & -a_1N_1(t)+2a_3N_3(t)(1-u_1(t))~, \\
   \dot{N}_2(t) & = & -a_2N_2(t)u_2(t)+a_1N_1(t)~, \\
   \dot{N}_3(t) & = & a_3N_3(t)+a_2N_2(t)u_2(t)~,
\end{array}
\end{array}
\end{equation}
the control inequality constraints
\begin{equation}\label{eq:threeCompartmentModel-ControlConstraint}
\begin{array}{rclcl}
  0 & \leq & u_1(t) & \leq & 1~, \\
  u_2^{\min} & \leq & u_2(t) & \leq & 1~, 
\end{array}
\end{equation}
and the boundary conditions
\begin{equation} \label{eq:threeCompartmentModel-BCs}
\begin{array}{lclclcl}
N_1(0) & = & 38~, &  & N_1(t_f) & = & \text{Free}~, \\
N_2(0) & = & 2.5~, & &  N_2(t_f) & = & \text{Free}~, \\
N_3(0) & = & 3.25 ~,&  & N_3(t_f) & = & \text{Free}~,
\end{array}
\end{equation} 
where $a_1=0.197$, $a_2=0.395$, $a_3=0.107$, $r_1=1$, $r_2=0.5$, $r_3=1$, $T=7$, and $u_2^{\min}=0.70$.  The optimal control problem given in Eqs.~\eqref{eq:threeCompartmentModel-Cost} -- \eqref{eq:threeCompartmentModel-BCs} was solved using each of the five mesh refinement methods $hp$-BB, $hp$-I, $hp$-II, $hp$-III, and $hp$-IV.  The solutions obtained using any of these five mesh refinement methods are in close agreement and match the solution given in Ref.~\cite{Ledzewicz1} (see pages $200$ -- $203$ of Ref.~\cite{Ledzewicz1}), and a summary of the performance of each method is shown in Table~\ref{tab:threeCompartmentModelData}.  In particular, it is seen in Table~\ref{tab:threeCompartmentModelData} that the $hp$-BB method is more computationally efficient, requires fewer mesh iterations, and results in a smaller final mesh to meet the mesh refinement accuracy tolerance of $10^{-6}$ when compared with any of the other mesh refinement methods.  Furthermore, as shown in Fig.~\ref{fig:threeCompartmentModelControlSolution}, the solution obtained using the $hp$-BB mesh refinement method accurately captures the bang-bang control profile of the optimal solution.  Specifically, it is seen that the optimal switch times are obtained to nearly machine precision.  Furthermore, estimates of the switching functions obtained using the solution on the initial mesh are shown in Fig.~\ref{fig:threeCompartmentModelSigma}, where it is seen that the roots of the switching functions are in close proximity to the locations of the optimal switch times.  Thus, while the previously developed $hp$-adaptive mesh refinement methods are able to satisfy the mesh refinement accuracy tolerance in a reasonable number of iterations, the $hp$-BB method outperforms all of these methods.  Moreover, unlike the other methods where a large number of collocation points are placed in the vicinity of a discontinuity (as seen in Fig.~\ref{fig:threeCompartmentModelControlSolution}), the $hp$-BB method places no unnecessary collocation points at the discontinuity due to the fact that optimal control problem has been divided into multiple domains and variable mesh points are included that define the locations of the discontinuities in the control (again, see Fig.~\ref{fig:threeCompartmentModelControlSolution}).

\begin{table}[!ht]
\centering
\renewcommand{\baselinestretch}{1}\small\normalfont
\caption{Mesh refinement performance results for Example 1 using $hp$-BB, $hp$-I, $hp$-II, $hp$-III, and $hp$-IV mesh refinement methods.
\label{tab:threeCompartmentModelData}}

\begin{tabular}{|c|c|c|c|c|c|}
\hline
	&	$hp$-BB	&	$hp$-I	&	~~~$hp$-II~~	&	~~~~$hp$-III~~~~	&	$hp$-IV		\\ \hline \hline
    $M$		&	$2$	&	$10$	&	$4$	&	$5$	&	$5$	\\ \hline
    $N_f$	&	$40$	&	$85$	&	$86$	&	$115$	&	$95$	\\ \hline
    $\C{T}$ (s)	&	$0.1234$	&	$0.4859$	&	$0.2826$	&	$0.4132$	&	$0.3054$	\\ \hline
\end{tabular}

\end{table}

\begin{figure}[!ht]
\centering

\includegraphics[keepaspectratio=true,width=2.5in]{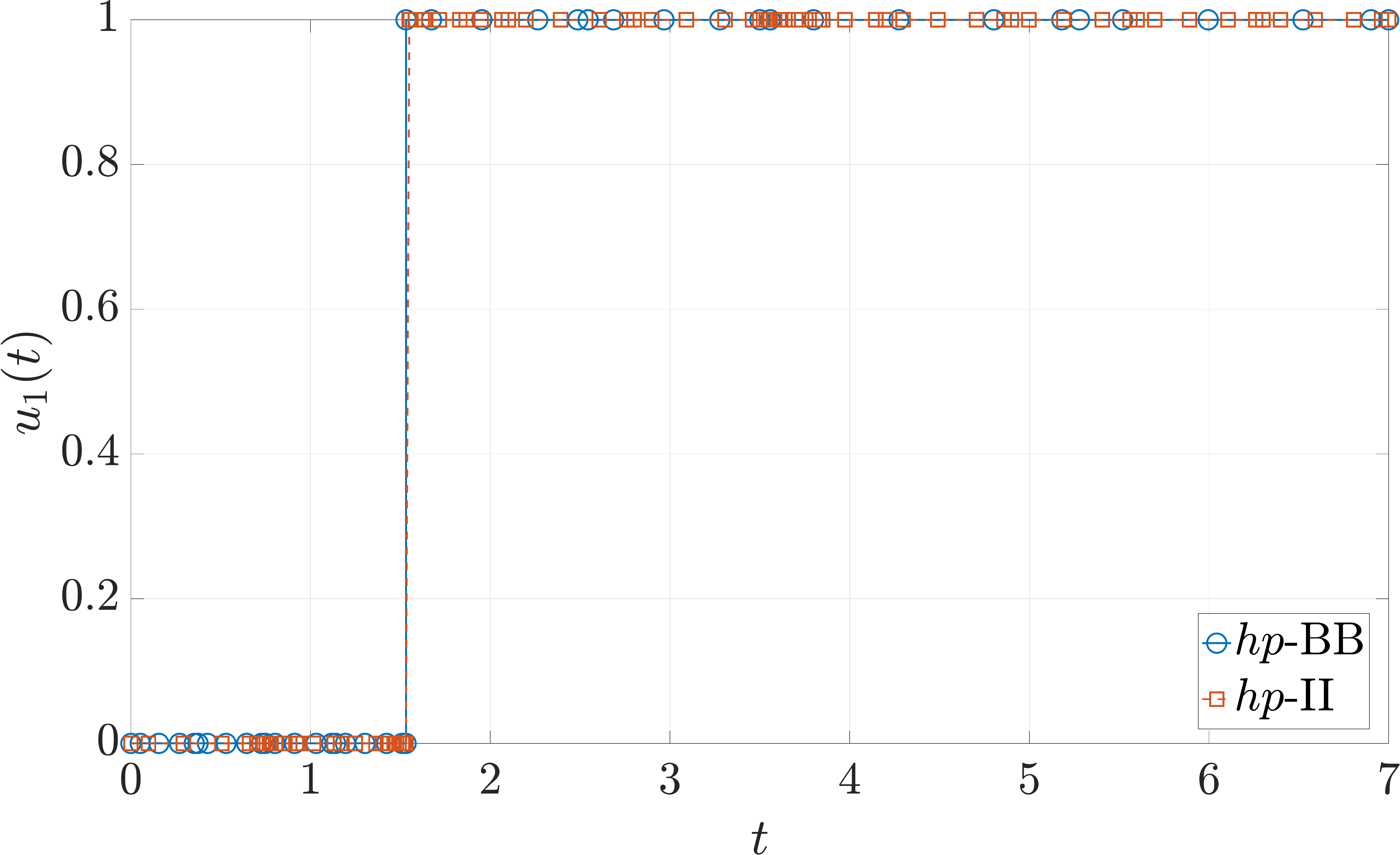}~~\includegraphics[keepaspectratio=true,width=2.5in]{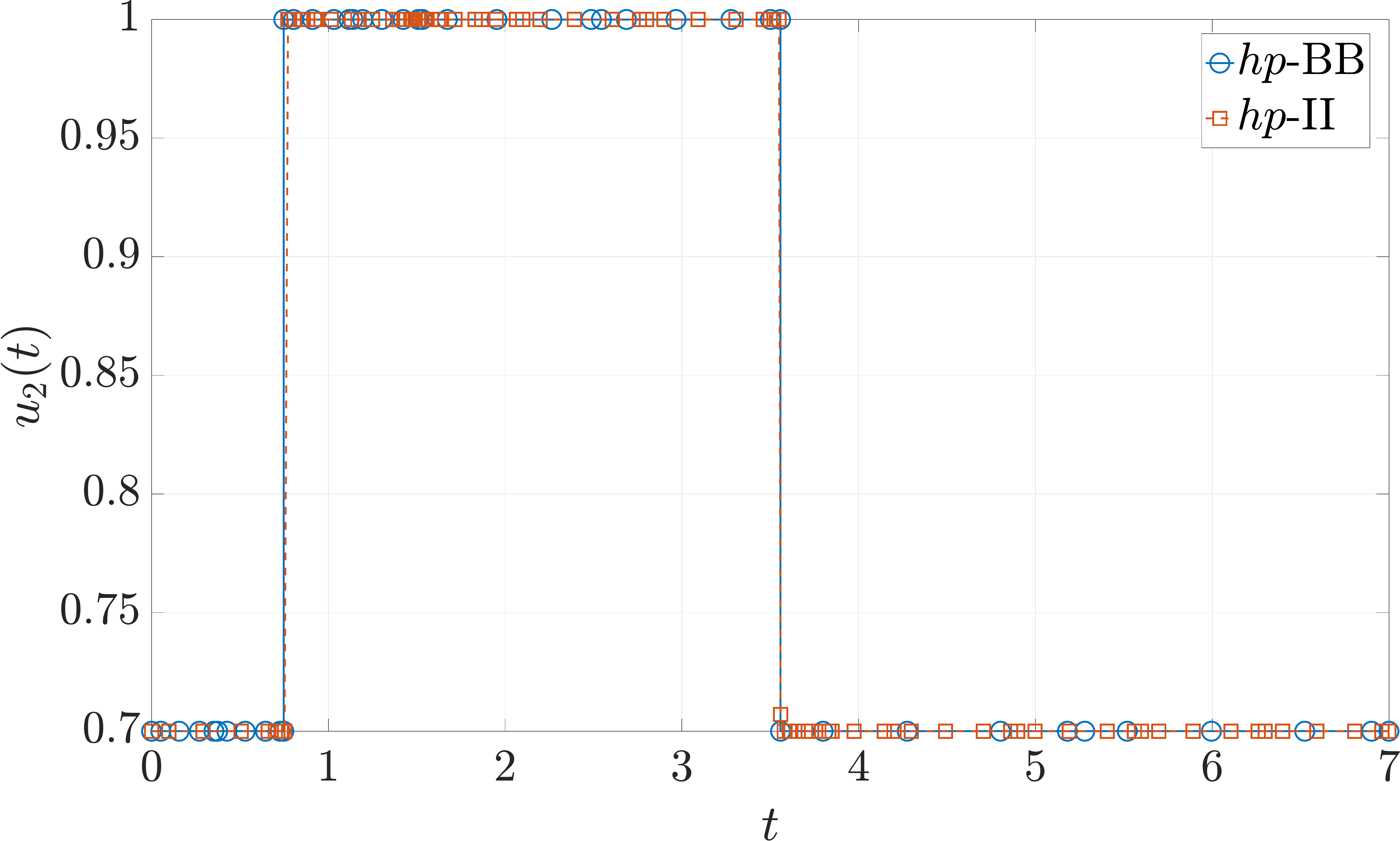}

\renewcommand{\baselinestretch}{1}\normalsize\normalfont
\caption{Comparison of control for Example 1 obtained using $hp$-BB and $hp$-II mesh refinement methods. \label{fig:threeCompartmentModelControlSolution}}

\end{figure}

\begin{figure}[!ht]
\centering

\includegraphics[keepaspectratio=true,width=2.5in]{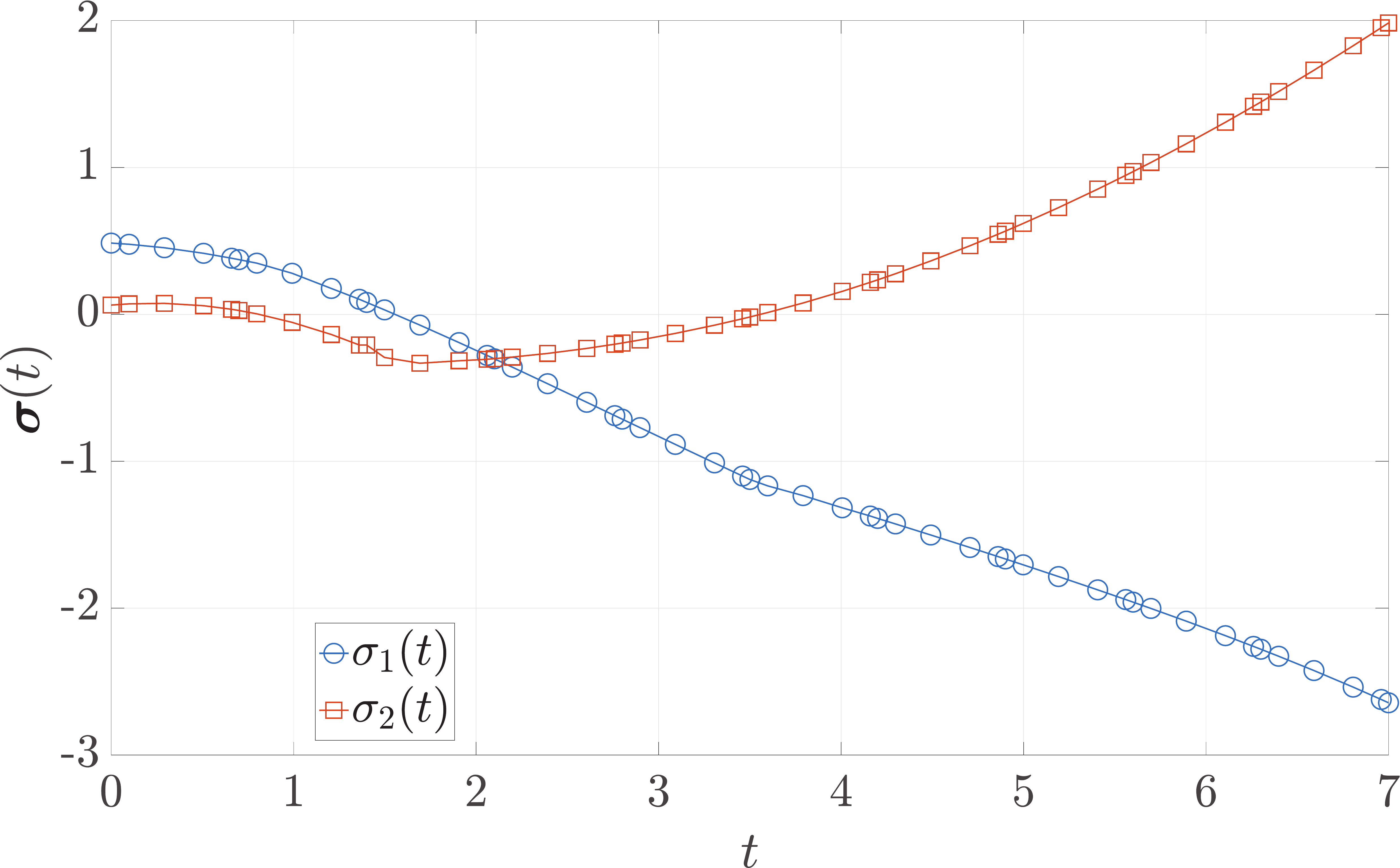}

\renewcommand{\baselinestretch}{1}\normalsize\normalfont
\caption{Estimates of the switching functions $\g{\sigma}(t)=(\sigma_1(t),\sigma_2(t))$ for Example 1 using solution obtained on the initial mesh. \label{fig:threeCompartmentModelSigma}}

\end{figure}

\subsection{Example 2: Robot Arm Problem \label{sect:robotArm}}

Consider the following optimal control problem taken from Ref.~\cite{Munson}.  Minimize the objective functional  
\begin{equation}\label{eq:robotArm-Cost}
  \C{J} = t_f~,
\end{equation}
subject to the dynamic constraints
\begin{equation}\label{eq:robotArm-Dynamics}
 \begin{array}{lclclclclcl}
   \dot{y}_1(t) & = & y_2(t)~, &  & \dot{y}_3(t) & = & y_4(t)~, & & \dot{y}_5(t) & = & y_6(t)~,\\
   \dot{y}_2(t) & = & u_1(t)/L~, &  & \dot{y}_4(t) & = & u_2(t)/I_{\theta}(t)~, & & \dot{y}_6(t) & = & u_3(t)/I_{\phi}(t)~,
\end{array}
\end{equation}
the control inequality constraints
\begin{equation}\label{eq:robotArm-ControlConstraint}
  -1 \leq u_i(t) \leq 1~, \quad (i=1,2,3)~, 
\end{equation}
and the boundary conditions
\begin{equation} \label{eq:robotArm-BCs}
\begin{array}{lclclclclclclcl}
 y_1(0) & = & 9/2~, &  & y_1(t_f) & = & 9/2~, & &
 y_2(0) & = & 0~, &  & y_2(t_f) & = & 0~, \\
 y_3(0) & = & 0 ~,&  & y_3(t_f) & = & 2\pi/3~, & &
 y_4(0) & = & 0~, &  & y_4(t_f) & = & 0~, \\
 y_5(0) & = & \pi/4~, &  & y_5(t_f) & = & \pi/4~, & &
 y_6(0) & = & 0~, &  & y_6(t_f) & = & 0~,
\end{array}
\end{equation} 
where 
\begin{equation}
\label{eq:robotArm-Aux1}
\begin{array}{ccc}
 I_{\theta}(t) =  \displaystyle \frac{((L - y_1(t))^3 + y_1^3(t))}{3}
 \sin^2(y_5(t))~, \quad& I_{\phi}(t)  =  \displaystyle \frac{((L - y_1(t))^3 +
   y_1^3(t))}{3}~, \quad&  L = 5~.
\end{array}
\end{equation}

The optimal control problem given in  Eqs.~\eqref{eq:robotArm-Cost}--\eqref{eq:robotArm-Aux1} was solved using each of the five mesh refinement methods $hp$-BB, $hp$-I, $hp$-II, $hp$-III, and $hp$-IV.  The solutions obtained using any of these five mesh refinement methods are in close agreement and match the solution given in Ref.~\cite{Munson} (see page $20$ of Ref.~\cite{Munson}), and a summary of the performance of each method is shown in Table~\ref{tab:robotArmData}.  In particular, it is seen in Table~\ref{tab:robotArmData} that the $hp$-BB method is more computationally efficient, requires fewer mesh iterations, and results in a smaller final mesh to meet the mesh refinement accuracy tolerance of $10^{-6}$ when compared with any of the other mesh refinement methods.  Furthermore, as shown in Fig.~\ref{fig:robotArmControlUSolution}, the solution obtained using the $hp$-BB mesh refinement method accurately captures the bang-bang control profile of the optimal solution.  Specifically, it is seen that the optimal switch times are obtained to nearly machine precision.  Furthermore, estimates of the switching functions obtained using the solution on the initial mesh are shown in Fig.~\ref{fig:robotArmSigma}, where it is seen that the roots of the switching functions are in close proximity to the locations of the optimal switch times.  Thus, while the previously developed $hp$-adaptive mesh refinement methods are able to satisfy the mesh refinement accuracy tolerance in a reasonable number of iterations, the $hp$-BB method outperforms all of these methods.  Moreover, unlike the other methods, where a large number of collocation points are placed in the vicinity of a discontinuity (as seen in Fig.~\ref{fig:robotArmControlUSolution}), the $hp$-BB method places no unnecessary collocation points at the discontinuity due to the fact that optimal control problem has been divided into multiple domains and variable mesh points are included that define the locations of the discontinuities in the control (again, see Fig.~\ref{fig:robotArmControlUSolution}).

\begin{table}[!ht]
\centering
\renewcommand{\baselinestretch}{1}\small\normalfont
\caption{Mesh refinement performance results for Example 2 using $hp$-BB, $hp$-I, $hp$-II, $hp$-III, and $hp$-IV mesh refinement methods.
\label{tab:robotArmData}}

\begin{tabular}{|c|c|c|c|c|c|}
\hline
	&	$hp$-BB	&	$hp$-I	&	~~~$hp$-II~~	&	~~~~$hp$-III~~~~	&	$hp$-IV		\\ \hline \hline
    $M$		&	$2$	&	$7$	&	$5$	&	$6$	&	$4$	\\ \hline
    $N_f$	&	$60$	&	$80$	&	$78$	&	$114$	&	$85$	\\ \hline
    $\C{T}$ (s)	&	$0.2183$	&	$0.6234$	&	$0.5000$	&	$0.8461$	&	$0.4156$	\\ \hline
\end{tabular}

\end{table}

\begin{figure}[!ht]
\centering

\includegraphics[keepaspectratio=true,width=2.5in]{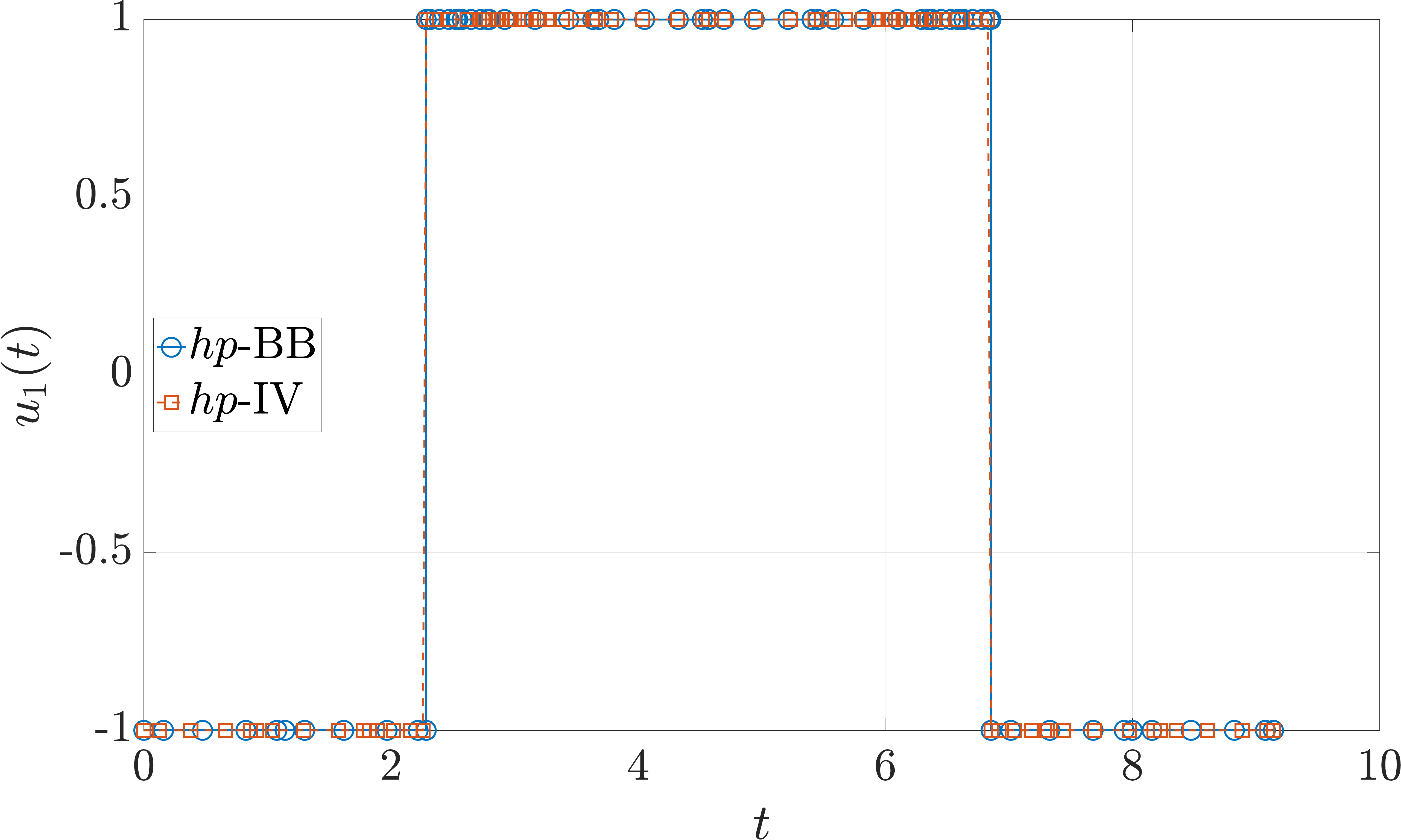}~~\includegraphics[keepaspectratio=true,width=2.5in]{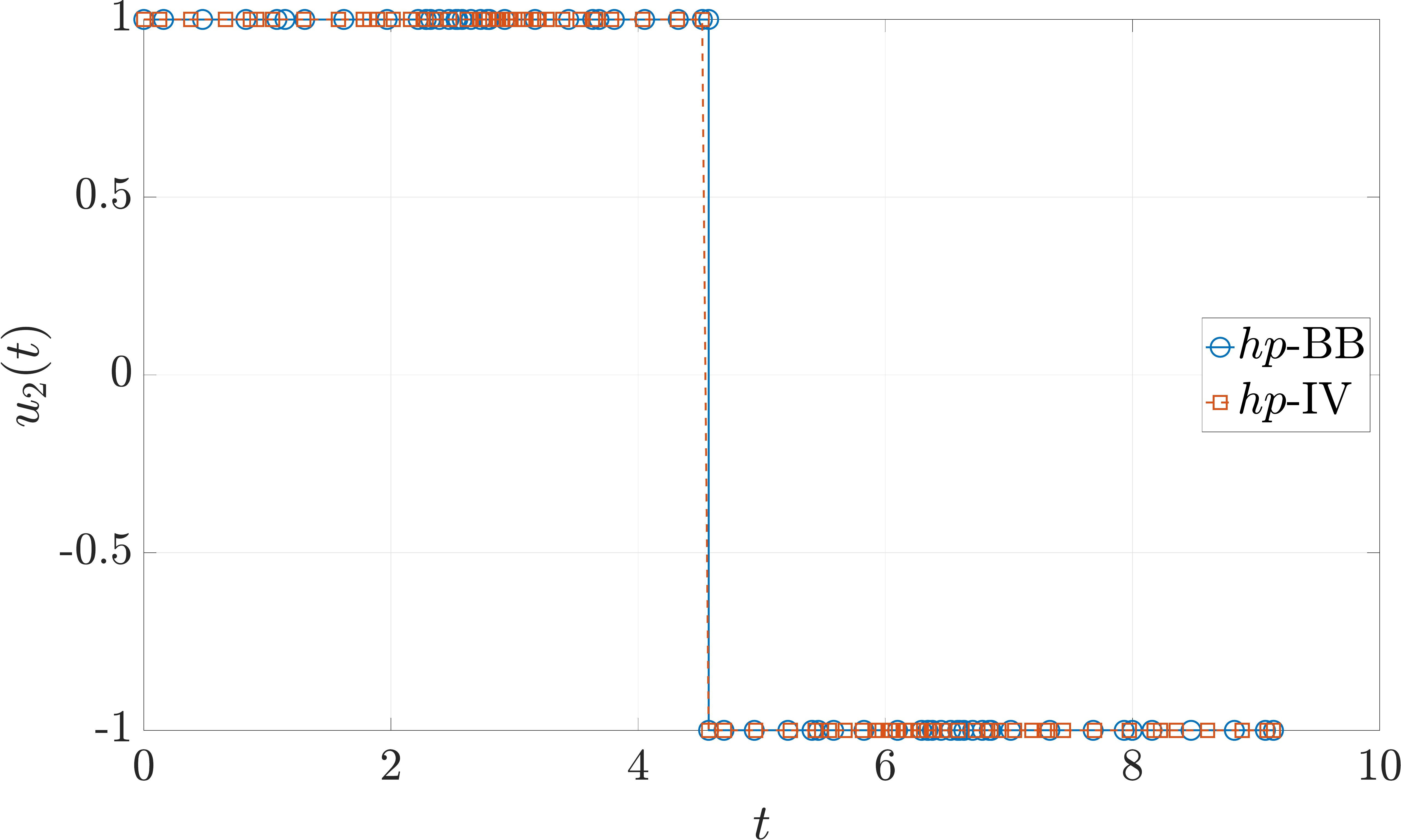}

\includegraphics[keepaspectratio=true,width=2.5in]{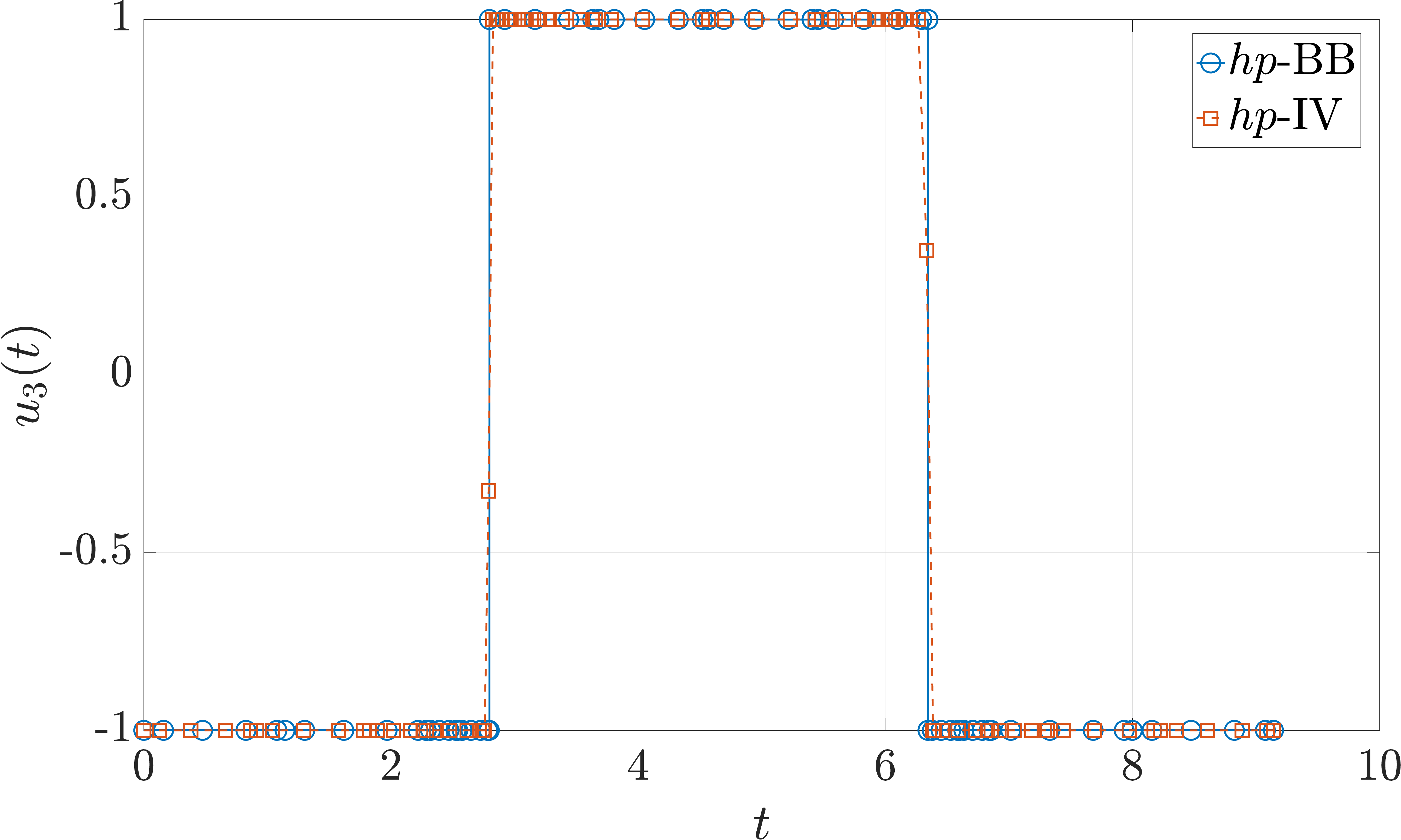}

\renewcommand{\baselinestretch}{1}\normalsize\normalfont
\caption{Comparison of control for Example 2 obtained Using $hp$-BB and $hp$-IV mesh refinement methods. \label{fig:robotArmControlUSolution}}

\end{figure}

\begin{figure}[!ht]
\centering

\includegraphics[keepaspectratio=true,width=2.5in]{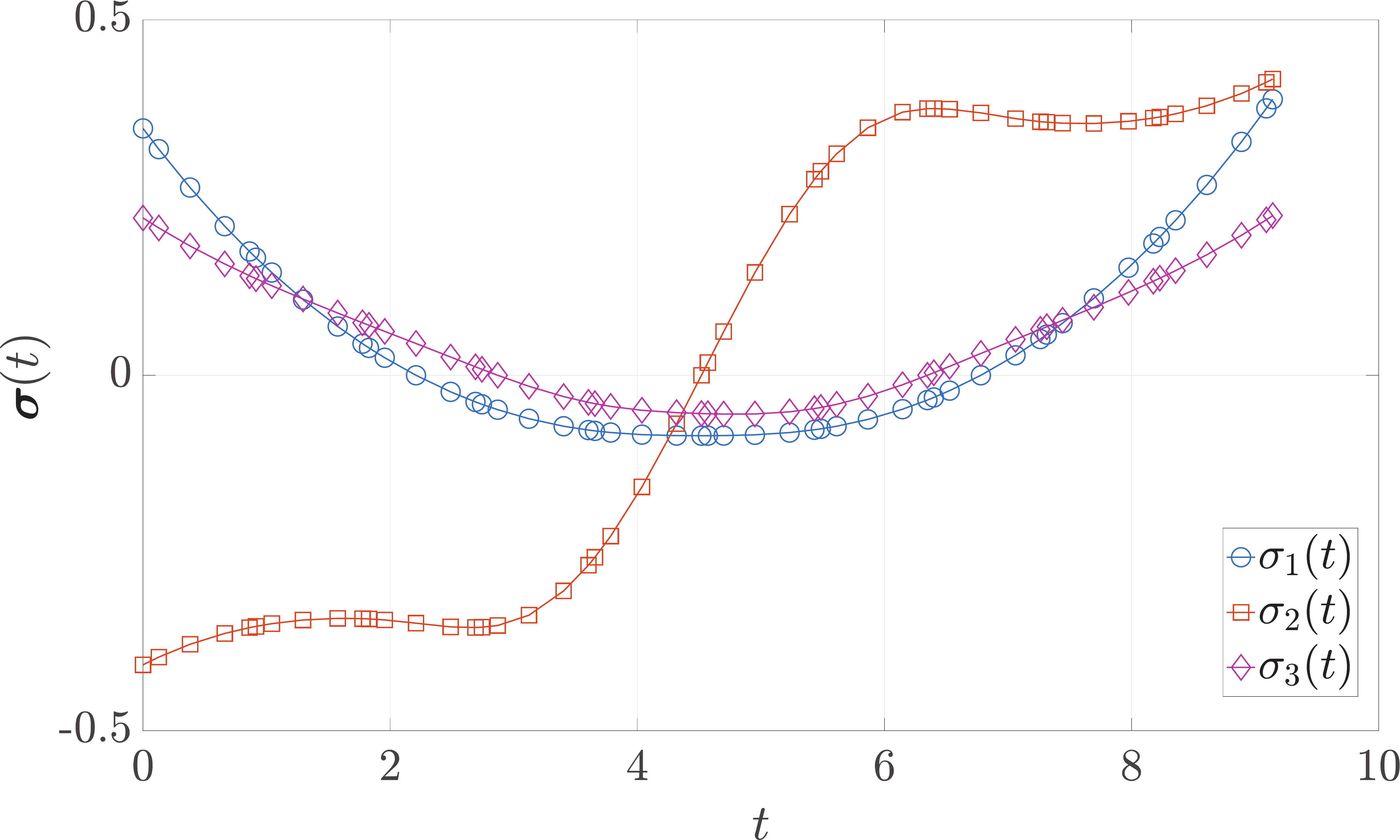}

\renewcommand{\baselinestretch}{1}\normalsize\normalfont
\caption{Estimates of the switching functions $\g{\sigma}(t)=(\sigma_1(t),\sigma_2(t),\sigma_3(t))$ for Example 2 using solution obtained on the initial mesh. \label{fig:robotArmSigma}}

\end{figure}

\clearpage

\subsection{Example 3: Free-Flying Robot Problem \label{sect:freeFlyingRobot}}

Consider the following optimal control problem taken from Ref.~\cite{Sakawa1}.  Minimize the cost functional  
\begin{equation}\label{eq:freeFlyingCost}
  \C{J} = \int_{0}^{t_f}(u_1(t)+u_2(t)+u_3(t)+u_4(t))dt~,
\end{equation}
subject to the dynamic constraints
\begin{equation}\label{eq:freeFlyingDynamics}
 \begin{array}{lclclcl}
   \dot{x}(t) & = & v_x(t)~, & & \dot{y}(t) & = & v_y(t)~, \\
   \dot{v_x}(t) & = & (F_1(t)+F_2(t))\cos(\theta(t))~, & & \dot{v_y}(t) & = & (F_1(t)+F_2(t))\sin(\theta(t))~, \\
   \dot{\theta}(t) & = & \omega(t)~, & & \dot{\omega}(t) & = & \alpha F_1(t)-\beta F_2(t)~, \\
\end{array}
\end{equation}
the control inequality constraints
\begin{equation}\label{eq:freeFlyingControlConstraint}
\begin{array}{lcl}
 0 \leq u_i(t) \leq 1~, \quad (i=1,2,3,4)~, & & \quad F_i(t) \leq 1~,~\quad (i=1,2)~,
\end{array}  
\end{equation}
and the boundary conditions
\begin{equation} \label{eq:freeFlyingBCs}
\begin{array}{lclclcl}
 x(0) = -10~, & & x(t_f) = 0~, & &
 y(0) = -10~, & & y(t_f) = 0~, \\
 v_x(0) = 0~, & & v_x(t_f) = 0~, & &
 v_y(0) = 0~, & & v_y(t_f) = 0~, \\
 \theta(0) = \frac{\pi}{2}~, & & \theta(t_f) = 0~, & &
 \omega(0) = 0~, & & \omega(t_f) = 0~,
\end{array}
\end{equation} 
where $F_1(t) = u_1(t)-u_2(t)$ and $F_2(t) = u_3(t)-u_4(t)$ are the real control, $\alpha = 0.2$, and $\beta = 0.2$.

The optimal control problem given in Eqs.~\eqref{eq:freeFlyingCost} -- \eqref{eq:freeFlyingBCs} was solved using each of the five mesh refinement methods $hp$-BB, $hp$-I, $hp$-II, $hp$-III, and $hp$-IV.  The solutions obtained using any of these five mesh refinement methods are in close agreement and match the solution given in Ref.~\cite{Munson} (see pages $328$ -- $329$ of Ref.~\cite{Betts3}), and a summary of the performance of each method is shown in Table~\ref{tab:freeFlyingRobotData}.  In particular, it is seen in Table~\ref{tab:freeFlyingRobotData} that the $hp$-BB method is more computationally efficient, requires fewer mesh iterations, and results in a smaller final mesh to meet the mesh refinement accuracy tolerance of $10^{-6}$ when compared with any of the other mesh refinement methods.  Furthermore, as shown in Fig.~\ref{fig:freeFlyingRobotRealControlSolution}, the solution obtained using the $hp$-BB mesh refinement method accurately captures the bang-bang control profile of the optimal solution.  Specifically, it is seen that the optimal switch times are obtained to nearly machine precision.  Furthermore, estimates of the switching functions obtained using the solution on the initial mesh are shown in Fig.~\ref{fig:freeFlyingRobotSigma}, where it is seen that the roots of the switching functions are in close proximity to the locations of the optimal switch times.  Thus, while the previously developed $hp$-adaptive mesh refinement methods are able to satisfy the mesh refinement accuracy tolerance in a reasonable number of iterations, the $hp$-BB method outperforms all of these methods.  Moreover, unlike the other methods, where a large number of collocation points are placed in the vicinity of a discontinuity (as seen in Fig.~\ref{fig:freeFlyingRobotRealControlSolution}), the $hp$-BB method places no unnecessary collocation points at the discontinuity due to the fact that optimal control problem has been divided into multiple domains and variable mesh points are included that define the locations of the discontinuities in the control (again, see Fig.~\ref{fig:freeFlyingRobotRealControlSolution}).

\begin{table}[!ht]
\centering
\renewcommand{\baselinestretch}{1}\small\normalfont
\caption{Mesh refinement performance results for Example 3 using $hp$-BB, $hp$-I, $hp$-II, $hp$-III, and $hp$-IV mesh refinement methods.
\label{tab:freeFlyingRobotData}}

\begin{tabular}{|c|c|c|c|c|c|}
\hline
	&	$hp$-BB	&	$hp$-I	&	~~~$hp$-II~~	&	~~~~$hp$-III~~~~	&	$hp$-IV		\\ \hline \hline
    $M$		&	$2$	&	$16$	&	$9$	&	$8$	&	$8$	\\ \hline
    $N_f$	&	$90$	&	$204$	&	$235$	&	$283$	&	$194$	\\ \hline
    $\C{T}$ (s)	&	$0.7019$	&	$8.4331$	&	$4.0833$	&	$4.7948$	&	$3.6157$	\\ \hline
\end{tabular}

\end{table}

\begin{figure}[!ht]
\centering

\includegraphics[keepaspectratio=true,width=2.5in]{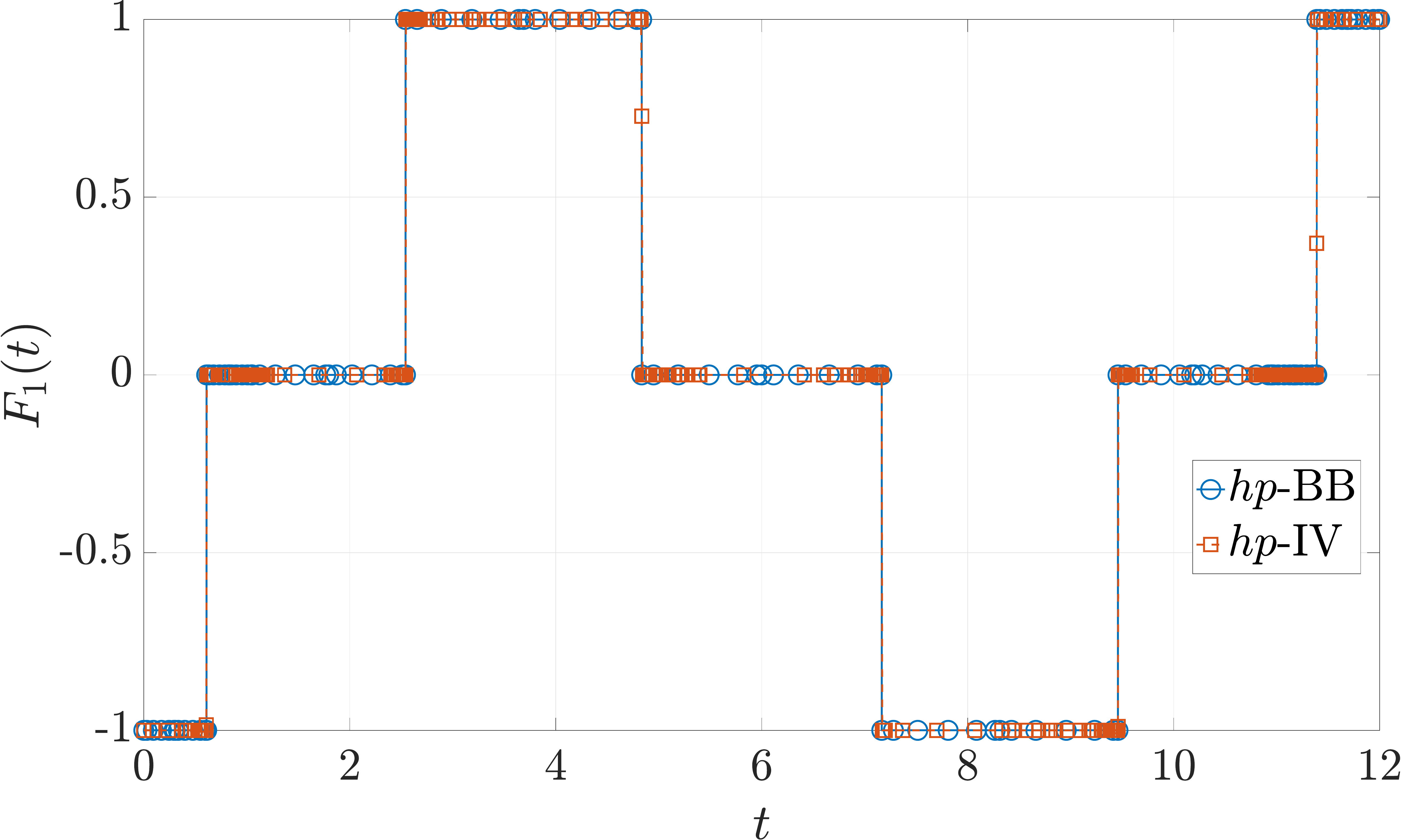}~~\includegraphics[keepaspectratio=true,width=2.5in]{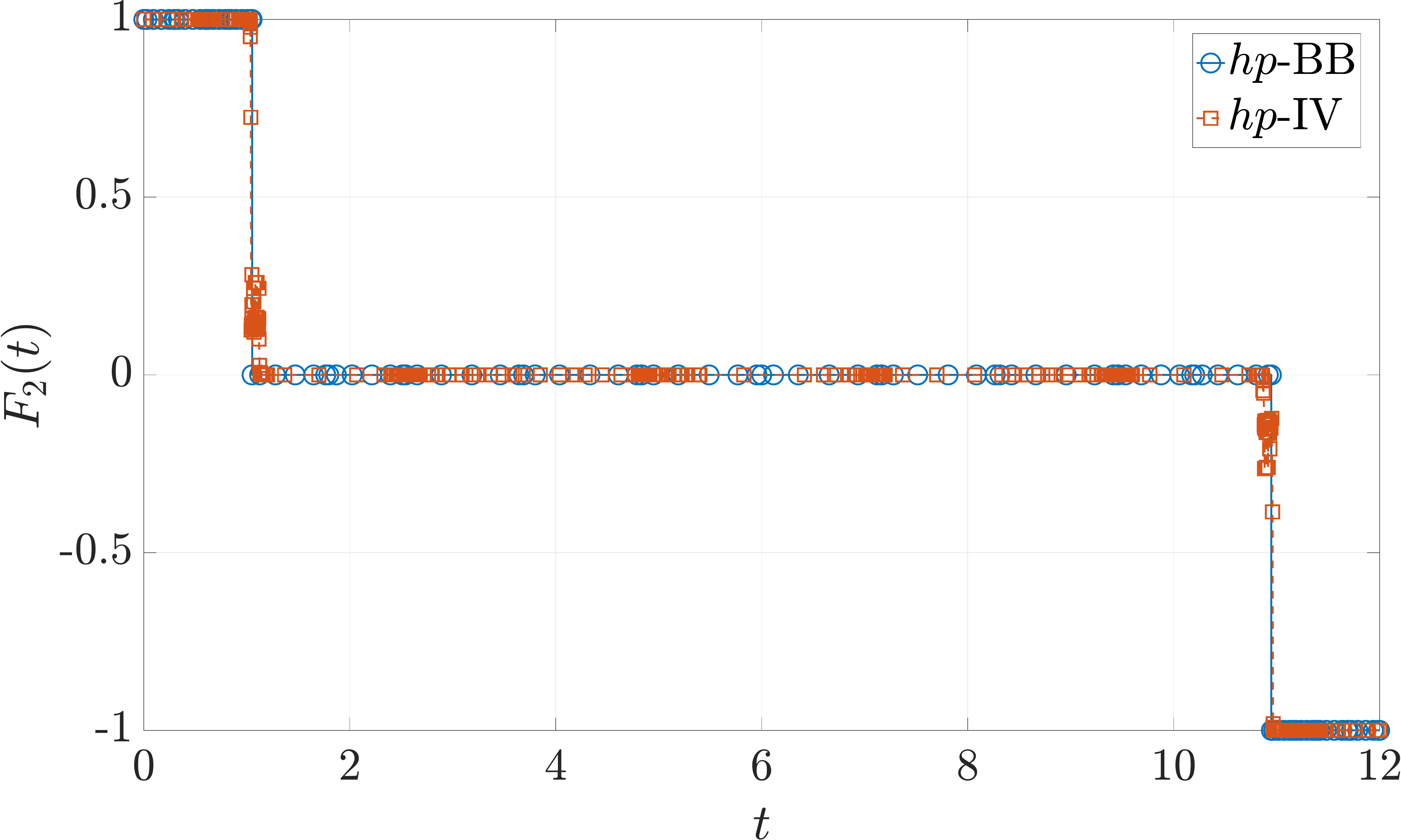}

\renewcommand{\baselinestretch}{1}\normalsize\normalfont
\caption{Comparison of control $(F_1(t),F_2(t))=(u_1(t)-u_2(t),u_3(t)-u_4(t))$ for Example 3 obtained using $hp$-BB, and $hp$-IV mesh refinement methods. \label{fig:freeFlyingRobotRealControlSolution}}

\end{figure}

\begin{figure}[!ht]
\centering

\includegraphics[keepaspectratio=true,width=2.5in]{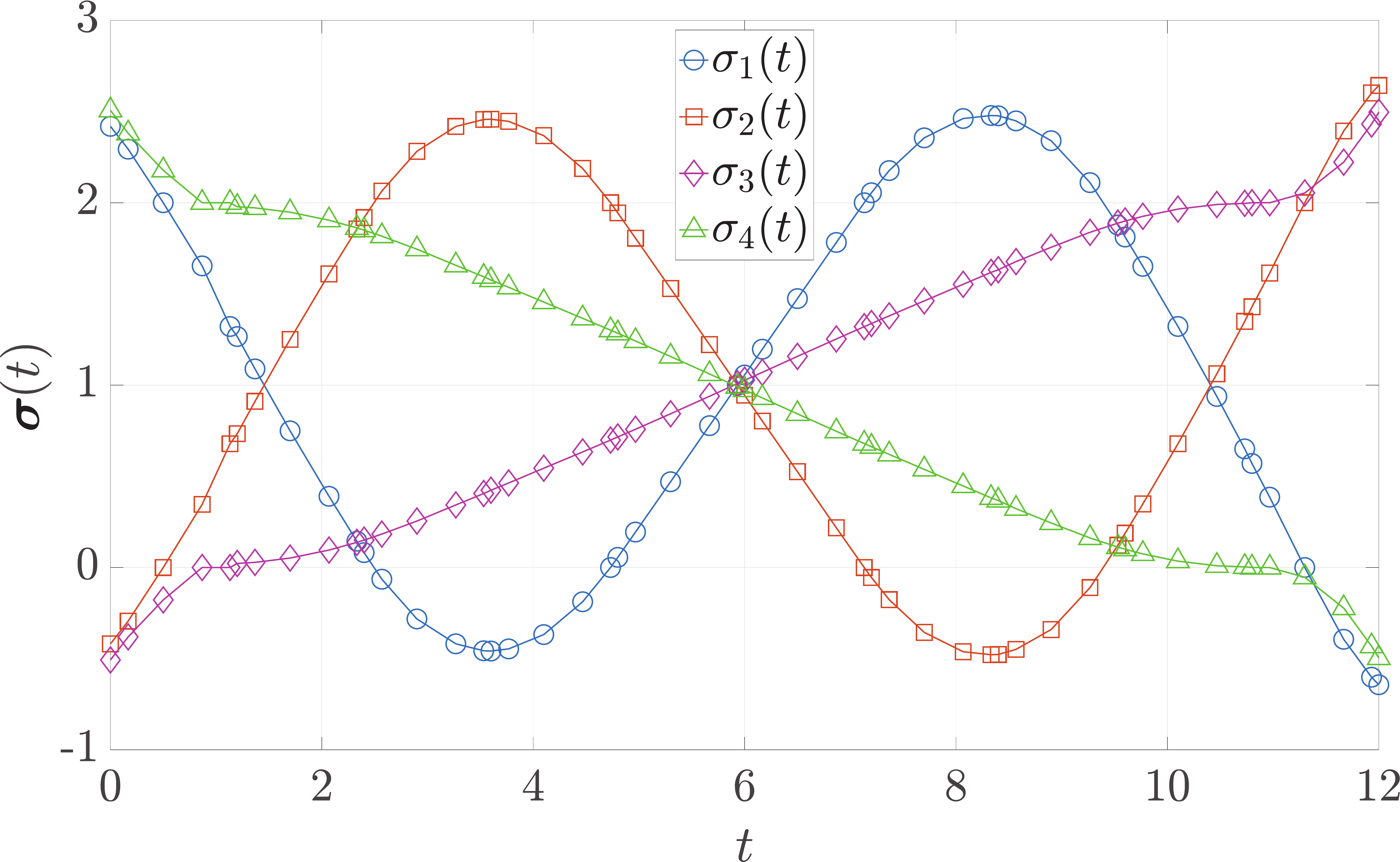}

\renewcommand{\baselinestretch}{1}\normalsize\normalfont
\caption{Estimates of the switching functions $\g{\sigma}(t)=(\sigma_1(t),\sigma_2(t),\sigma_3(t),\sigma_4(t))$ for Example 3 using solution obtained on the initial mesh. \label{fig:freeFlyingRobotSigma}}

\end{figure}

\clearpage

\section{Discussion \label{sect:Discussion}}

The results of Section \ref{sect:Examples} demonstrate the effectiveness of the mesh refinement method developed in this paper for problems whose optimal control has a bang-bang structure.  In particular, the results of Section \ref{sect:Examples} show that, while the previously developed mesh refinement methods are able to find a solution that meets a specified mesh refinement accuracy tolerance, these methods place an unnecessarily large number of collocation points in the vicinity of a discontinuity in the control.  In addition, these methods often require a large amount of mesh refinement to meet a desired accuracy tolerance.  On the other hand, for problems whose optimal control has a bang-bang structure, the mesh refinement method developed in this paper locates discontinuities accurately.  This improved accuracy (over a standard mesh refinement method) is due to the fact that accurate estimates are obtained of the switching functions associated with those components of the control upon which the Hamiltonian depends linearly (where an accurate estimate of the switching functions is obtained because the costate of the optimal control problem is approximated accurately using the LGR collocation method).  Then, by partitioning the horizon into multiple domains, introducing variables that define the locations of the switch times in the control-linear components, and fixing the control-linear components to lie at either its lower or upper limit in each domain, the method developed in this paper accurately identifies the switch times.   Moreover, solutions that meet the specified accuracy tolerance are obtained in fewer mesh refinement iterations when compared with using one of the previously developed mesh refinement methods.  

\section{Conclusions \label{sect:Conclusions}}

A mesh refinement method has been described for solving bang-bang optimal control problems using direct collocation.  First, the solution of the optimal control problem is approximated on a coarse mesh.  If the approximation on the coarse mesh does not satisfy the specified accuracy tolerance, the method determines automatically if the Hamiltonian of the optimal control problem is linear with respect to any components of the control.  Then, for any control component upon which the Hamiltonian depends linearly, the locations of the discontinuities in the control are obtained by estimating the roots of the switching functions associated with any component of the control that appears linearly in the Hamiltonian using estimates of the switching functions obtained using the state and costate obtained from the solution on the initial mesh.  The estimates of the switching functions are then used to determine the bang-bang structure of the optimal solution.   Using estimates of the locations of discontinuities in the control obtained from the detected structure, the horizon is partitioned into multiple domains and parameters corresponding to the the locations of the switch times are introduced as variables in the optimization.  Then, by fixing in each domain any component of the control that has a bang-bang structure to lie at either its lower or upper limit, the multiple-domain optimal control problem is solved to accurately determine the switch times.  The method has been demonstrated on three examples where it has been shown to efficiently obtain accurate approximations to solutions of bang-bang optimal control problems.

\section*{Acknowledgments}

The authors gratefully acknowledge support for this research from the U.S.~Office of Naval Research under grant N00014-15-1-2048 and from the U.S.~National Science Foundation under grants DMS-1522629, DMS-1924762, and CMMI-1563225.

\renewcommand{\baselinestretch}{1.0}
\normalsize\normalfont
\bibliographystyle{aiaa}

\end{document}